\documentclass[11pt,a4paper]{amsart}

\theoremstyle{plain}
\newtheorem{thm}{Theorem}[section]
\newtheorem*{thm*}{Theorem}
\newtheorem{lem}[thm]{Lemma}
\newtheorem{prop}[thm]{Proposition}
\newtheorem{cor}[thm]{Corollary}
\theoremstyle{definition}

\newtheorem{exmp}[thm]{Example}

\theoremstyle{remark}
\newtheorem{rem}[thm]{Remark}

\usepackage{amsfonts, amsmath, amssymb, amsthm}
\usepackage{amsmath, amssymb}
\usepackage{bbm}

\usepackage{color}
\usepackage{graphicx}
\usepackage{url}
\usepackage{overpic}

\usepackage{psfrag}
\usepackage{paralist}
\usepackage[a4paper]{anysize}
\usepackage[all,ps]{xy}

\newcommand\subgr{\operatorname{<}}

\newcommand\Sing{\operatorname{sing}}
\newcommand{\Sym}{\operatorname{Sym}}

\newcommand{\s}{\sigma}
\newcommand{\g}{\gamma}
\newcommand{\Star}{\operatorname{st}}
\newcommand{\Link}{\operatorname{lk}}
\newcommand{\Cone}{\operatorname{cone}}
\newcommand{\Susp}{\operatorname{susp}}

\newcommand{\Codim}{\operatorname{codim}}

\newcommand{\mysim}{\hspace{-.14cm}\sim}
\newcommand{\B}{\mathbbm{D}}
\newcommand{\Id}{\operatorname{Id}}
\newcommand{\M}{\mathfrak{M}}
\newcommand{\m}{\mathfrak{m}}
\newcommand{\h}{\mathfrak{h}}
\newcommand\SetOf[2]{\left\{#1\,\vphantom{#2}\right.\left|\vphantom{#1}\,#2\right\}}

\newcommand{\odd}{\operatorname{odd}}
\newcommand{\closure}{\operatorname{cl}}
\newcommand{\dotcup}{\ensuremath{\mathaccent\cdot\cup}}
\newcommand\RR{{\mathbb R}}
\newcommand\CC{{\mathbb C}}
\newcommand\Sph{{\mathbb S}}
\newcommand\restr{\kern-1ex\mid}

\begin{document}

\title{Constructing Combinatorial 4-Manifolds}
\author[Witte]{Nikolaus Witte}
 \address{Nikolaus Witte, Fachbereich Mathematik, AG~7, TU Darmstadt, 64289 Darmstadt, Germany}
 \email{witte@math.tu-berlin.de}
 \thanks{The authors is supported by Deutsche Forschungsgemeinschaft, DFG Research Group
   ``Polyhedral Surfaces''.}
 \date{\today}

\begin{abstract}
  Every closed oriented PL 4-manifold is a branched cover of the 4-sphere branched over a
  PL-surface with finitely many singularities by Piergallini [Topology 34(3):497-508, 1995].
  This generalizes a long standing result by Hilden and Montesinos to dimension four. Izmestiev and
  Joswig [Adv.\ Geom.\ 3(2):191-225, 2003] gave a combinatorial equivalent of the Hilden and
  Montesinos result, constructing closed oriented combinatorial 3-manifolds as simplicial branched
  covers of combinatorial 3-spheres. The construction of Izmestiev and Joswig is generalized and
  applied to the result of Piergallini, obtaining closed oriented combinatorial 4-manifolds as
  simplicial branched covers of simplicial 4-spheres.
\end{abstract}


\maketitle


\section{Introduction}

\noindent
The main objective of this paper is to give a complete yet concise
account of how to obtain closed oriented combinatorial 4-manifolds as
simplicial branched covers, that is, as partial unfolding of
simplicial 4-spheres. The construction is at times technical involved
and the topological background extensive. Thus we abstain from
discussing related material and omit some of the proofs. Complete
proofs and plenty of further material can be found in the Chapters~1
to~3 of~\cite{witte:DISS}.  Concerning the construction of closed
oriented PL 4-manifolds as branched covers we refer to
Piergallini~\cite{piergallini:FMS4} and
Montesinos~\cite{montesinos:CSR}.  For the partial unfolding and the
construction of closed oriented combinatorial 3-manifolds Izmestiev~\&
Joswig~\cite{izmestiev_joswig:BC} is mandatory reading. (Their
construction has recently be simplified significantly by Hilden,
Montesinos-Amilibia, Tejada~\& Toro~\cite{MR2218370}.) For those able
to read German additional analysis and examples can be found
in~\cite{witte:DIPL}. The partial unfolding is implemented in the
software package \texttt{polymake}~\cite{polymake}.

Branched covers form a major tool for the study, construction and
classification of $d$-manifolds. First results are by
Alexander~\cite{alexander:NORS} in 1920, who observed
that any closed
oriented PL $d$-manifold~$M$ is a branched cover
of the $d$-sphere.
%
Unfortunately Alexander's proof does not allow for any (reasonable) control
over the number of sheets of the branched cover, nor over the topology
of the branching set:
The number of sheets depends on the size of some
triangulation of~$M$ and the branching set is the co-dimension 2-skeleton of the $d$-simplex.

At least to our knowledge, there are no
non-trivial upper bounds for the number of sheets of such a branched
cover for $d>4$. On the contrary, Bernstein~\&
Edmonds~\cite{bernstein_edmonds:DABS} showed that at least~$d$ sheets
are necessary in general (for example the $d$-torus $(\Sph^1)^d$
exhibits such a behavior), and that the branching set can not be
required to be non-singular for $d \geq 8$.

However, in dimension $d\leq 4$, the situation is fairly well
understood. The 2-dimensional case is straight forward;
any closed oriented surface~$F_g$ of genus~$g$ is a
2-fold branched cover of the 2-sphere branched over $2g+2$
isolated points.

By results of Hilden~\cite{hilden:TFBC} and
Montesinos~\cite{montesinos:MFBC} any closed oriented 3-manifold~$M$
arises as 3-fold simple branched cover of the 3-sphere branched over a
link~$L$. Labeling each bridge~$b$ of a diagram of~$L$ with the
corresponding monodromy action of a meridian around~$b$, we can
represent~$M$ as a labeled (or colored) link diagram.


In dimension~4 the situation becomes increasingly difficult.
First Piergallini~\cite{piergallini:FMS4}
showed how to obtain any closed oriented 
PL 4-manifold as a 4-fold branched cover
of the 4-sphere branched over a PL-surface
with a finite number of cusp and node singularities.
Prior to
Piergallini's work Montesinos~\cite{montesinos:CSR} gave a description
of oriented 4-manifolds composed \mbox{of~0-,~1-,} and 2-handles only as a
branched cover of the 4-ball.  Montesinos' result is essential for
Piergallini's construction of closed oriented PL 4-manifolds as
branched covers.
These two constructions are the ``blue
print'' for the main result of this paper and they are reviewed in
Section~\ref{sec:4mf_as_branched_covers}.

Piergallini and later Iori~\& Piergallini improved the results on the
construction of closed oriented PL 4-manifolds further. First
Piergallini~\cite{piergallini:FMS4} eliminated the cusp singularities
of the branching set.  This yields a branched cover with a
transversally immersed PL-surface as its branching set. Iori~\&
Piergallini~\cite{iori_piergallini:4MFNS} then proved that the
branching set may be realized locally flat if one allows for a fifth
sheet for the branched cover, thus proving a long-standing conjecture
by Montesinos~\cite{montesinos:CSR}. The question whether any closed
oriented PL 4-manifold can be obtained as 4-fold cover of the 4-sphere
branched over a locally flat PL-surface is still open. Although these
later developments certainly ask for a combinatorial equivalent, we
will not investigate these here, nor make use of these observations.

\subsubsection*{Outline of the paper.}
After some basic definitions and notations the partial
unfolding~$\widehat{K}$ of a simplicial complex~$K$ is introduced. The
partial unfolding defines a projection $p:\widehat{K}\to K$ which is a
simplicial branched cover if~$K$ meets certain connectivity
assumptions. We define combinatorial models of key features of a
branched cover, namely the branching set and the monodromy
homomorphism.

Section~\ref{sec:color_equivalence} introduces a notion of equivalence of simplicial
complexes which agrees with their unfolding behavior. We proceed by
establishing further (technical) tools for the construction of
combinatorial 4-manifolds in Section~\ref{sec:constructing_4mfs}.

Finally Section~\ref{sec:constructing_4mfs} states and proofs the main result
Theorem~\ref{thm:ch_4mf_main}.
The key idea is to construct a simplicial 4-sphere~$S$, such
that the projection $p:\widehat{S}\to S$ is equivalent to
a given branched cover $r:M\to\Sph^4$. In particular, the equivalence of the
branched covers~$p$ and~$r$ implies homeomorphy of the covering
spaces~$\widehat{K}$ and~$M$.
In Theorem~\ref{thm:ch_4mf_main} we prove that this is indeed
possible for the branched covers arising in the construction of closed
oriented PL 4-manifolds by Piergallini~\cite{piergallini:FMS4}: For any given closed oriented PL 4-manifold~$M$ there is
a simplicial 4-sphere~$S$ such that the partial
unfolding~$\widehat{S}$ is PL-homeomorphic to~$M$.  We proceed by
giving a construction of the simplicial 4-sphere~$S$.
Prior to proving the main result, the topological
constructions by Montesinos~\cite{montesinos:CSR} and Piergallini~\cite{piergallini:FMS4} are
reviewed.

\subsection{Basic definitions and notations.}

Given some topological manifold~$M$, we call a simplicial complex~$K$
homeomorphic to~$M$ a \emph{triangulation} of~$M$, or a
\emph{simplicial manifold}. A simplicial complex~$K$ is a
\emph{combinatorial $d$-sphere} or \emph{combinatorial $d$-ball} if it
is piecewise linear homeomorphic to the boundary of the
$(d+1)$-simplex, respectively to the $d$-simplex. Equivalently,~$K$ is
a combinatorial $d$-sphere or $d$-ball if there is a common refinement
of~$K$ and the boundary of the $(d+1)$-simplex, respectively the
$d$-simplex. A simplicial complex~$K$ is a \emph{combinatorial
  manifold} if the vertex link of each vertex of~$K$ is a
combinatorial sphere or a combinatorial ball. Note that combinatorial
spheres and balls are combinatorial manifolds.

A manifold~$M$ where all charts are piecewise linear is called a
\emph{PL-manifold}. Up to dimension~3 there is no difference between
topological, PL-, and differential manifolds, that is, every
topological manifold allows for a PL- or differential atlas (or
structure). The existence of a triangulation of~$M$ as a combinatorial
manifold is equivalent to the existence of a PL-atlas for~$M$.  For an
introduction to PL-topology see Bj\"orner~\cite[Part~II]{bjoerner:TM},
Hudson~\cite{hudson:PLT}, and Rourke~\&
Sanderson~\cite{rourke_sanderson:PLT}.

Similarly to the topological situation, there is no difference between
the notion of a simplicial and a combinatorial manifold in
dimension~$d\leq 3$, that is, every simplicial manifold (or sphere, or
ball) is a combinatorial manifold (or sphere, or ball). But in
dimension~4 the situation becomes more complicated. Freedman~\&
Quinn~\cite{freedman_quinn:TO4MF} construct a 4-manifold which does
not have a triangulation as a combinatorial manifold. In fact, there
are 4-manifolds which can not be triangulated at
all~\cite[p.~9]{lutz:MFV}. The following unanswered question
illustrates the subtleties of the 4-dimensional case like no other: Is
a combinatorial manifold homeomorphic to the 4-sphere necessarily a
combinatorial 4-sphere? Surprisingly, the answer to this question is
affirmative in all dimensions $d\not=4$; see
Moise~\cite{moise:AS3MF-V} and Kirby~\&
Siebenmann~\cite{kirby_siebenmann:ETMF}.

Neither barycentric subdivision nor anti-prismatic subdivision (of a
face) change the PL-type of a simplicial manifold, that is, the
subdivision of a simplicial complex~$K$ is a combinatorial manifold if
and only if~$K$ is a combinatorial manifold. The cone of a
combinatorial sphere is a combinatorial ball and the suspension of a
combinatorial sphere is again a combinatorial sphere.

The simplicial complexes considered in the following (and throughout
this exposition) are always \emph{pure}, that is, all the inclusion
maximal faces, called the \emph{facets}, have the same dimension.  We
call a co-dimension 1-face of a pure simplicial complex a
\emph{ridge}, and the \emph{dual graph}~$\Gamma^*(K)$ of a pure
simplicial complex~$K$ has the facets as its node set, and two nodes
are adjacent if the corresponding facets share a ridge. We denote the
1-skeleton of~$K$ by~$\Gamma (K)$, its \emph{graph}.

Further it is often necessary to restrict ourselves to simplicial
complexes with certain connectivity properties: A pure simplicial
complex~$K$ is \emph{strongly connected} if its dual
graph~$\Gamma^*(K)$ is connected, and \emph{locally strongly
  connected} if the star $\Star_K(f)$ of~$f$ is strongly connected for
each face $f\in K$. If~$K$ is locally strongly connected, then
connected and strongly connected coincide. Further we call~$K$
\emph{locally strongly simply connected} if for each face $f\in K$
with co-dimension~$\geq2$ the link~$\Link_K(f)$ of~$f$ is simply
connected, and finally,~$K$ is \emph{nice} if it is locally strongly
connected and locally strongly simply connected.  Observe that
connected combinatorial manifolds are always nice.

\subsection{The branched cover.}
\label{sec:bc}

The concept of a covering of a space~$Y$ by another space~$X$ is
generalized by Fox~\cite{fox:CSWS} to the notion of the branched
cover. Here a certain subset~$Y_{\Sing}\subset Y$ may violate the
conditions of a covering map. This allows for a wider application in
the construction of topological spaces.  It is essential for a
satisfactory theory of (branched) coverings to make certain
connectivity assumption for~$X$ and~$Y$. The spaces mostly considered
are Hausdorff, path connected, and locally path connected; see
Bredon~\cite[III.3.1]{bredon:TAG}.  Throughout we will restrict our
attention to coverings of manifolds and we assume~$Y$ to be
connected, hence they meet the connectivity assumptions
in~\cite{bredon:TAG}.

Consider a continuous map $h:Z\to Y$, and assume the
restriction~$h:Z\to h(Z)$ to be a covering. If~$h(Z)$ is dense in~$Y$
(and meets certain additional connectivity conditions) then there is a
surjective map $p:X\to Y$ with $Z\subset X$ and
$p\left|_{Z}\right.=h$. The map~$p$ is called a \emph{completion}
of~$h$, and any two completions $p:X\to Y$ and $p':X'\to Y$ are
equivalent in the sense that there exists a homeomorphism
$\varphi:X\to X'$ satisfying $p'\circ\varphi=p$ and
$\varphi\left|_Z\right.=\Id$.  The map~$p:X\to Y$ obtained this way is
a \emph{branched cover}, and we call the unique minimal subset
$Y_{\Sing}\subset Y$ such that the restriction of~$p$ to the preimage
of $Y\setminus Y_{\Sing}$ is a cover, the \emph{branching set} of~$p$.
The restriction of~$p$ to $p^{-1}(Y\setminus Y_{\Sing})$ is called the
\emph{associated} cover of~$p$. If $h:Z\to Y$ is a cover, then~$X=Z$,
and~$p=h$ is a branched cover with empty branching set.

\begin{exmp}\label{exmp:branched_cover}
  For $k\geq 2$ consider the map $p_k:\CC \to \CC$. The restriction
  $p_k\left|_{\B^2}\right.$ is a $k$-fold branched cover $\B^2 \to
  \B^2$ with the single branch point~$\{0\}$.
\end{exmp}

We define the \emph{monodromy homomorphism}
\[
\m_p:\pi_1(Y\setminus Y_{\Sing},y_0)\to\Sym(p^{-1}(y_0))
\]
of a branched cover for a point $y_0\in Y\setminus Y_{\Sing}$ as the
monodromy homomorphism of the associated cover: If
$[\alpha]\in\pi_1(Y\setminus Y_{\Sing},y_0)$ is represented by a closed
path~$\alpha$ based at~$y_0$, then~$\m_p$ maps~$[\alpha]$ to the
permutation $(x_i\mapsto \alpha_i(1))$, where $\{x_1,x_2,\dots,x_k\}=p^{-1}(y_0)$
is the preimage of~$y_0$ and $\alpha_i:[0,1]\to X$ is the unique
lifting of~$\alpha$ with $p\circ\alpha_i=\alpha$ and $\alpha_i(0)=x_i$; see
Munkres~\cite[Lemma~79.1]{munkres:TOP} and Seifert~\& Threlfall~\cite[\S~58]{seifert_threlfall:TOT}.
The \emph{monodromy group}~$\M_p$ is defined as the image of~$\m_p$.

Two branched covers $p:X\to Y$ and~$p':X'\to Y'$ are \emph{equivalent}
if there are homeomorphisms $\varphi:X\to X'$ and $\psi:Y\to Y'$ with
$\psi(Y_{\Sing})=Y'_{\Sing}$, such that $p'\circ\varphi=\psi\circ p$
holds.  The well known Theorem~\ref{thm:bc} is due to the uniqueness
of~$Y_{\Sing}$, and hence the uniqueness of the associated cover;
see Piergallini~\cite[p.~2]{piergallini:MBC}.

\begin{thm}\label{thm:bc}
  Let $p:X\to Y$ be a branched cover of a connected manifold~$Y$.
  Then~$p$ is uniquely determined up to equivalence
  by the branching set~$Y_{\Sing}$, and the monodromy
  homeomorphism~$\m_p$. In particular, the covering space~$X$ is
  determined up to homeomorphy.
\end{thm}

Let~$Y$ be a connected manifold and~$Y_{\Sing}$ a co-dimension~2
submanifold, possibly with a finite number of singularities. We call a
branched cover~$p$ \emph{simple} if the image~$\m_p(m)$ of any
meridial loop~$m$ around a non-singular point of the branching set is
a transposition in~$\M_p$. Note that the $k$-fold branched cover
$p_k\left|_{\B^2}\right. :\B^2\to\B^2$ presented in
Example~\ref{exmp:branched_cover} is not simple for $k\geq3$.

\subsection{The partial unfolding.}

The partial unfolding~$\widehat{K}$ of a simplicial complex~$K$ first
appeared in a paper by Izmestiev~\& Joswig~\cite{izmestiev_joswig:BC},
with some of the basic notions already developed in
Joswig~\cite{MR1900311}. The partial unfolding is closely related to
the complete unfolding, also defined in~\cite{izmestiev_joswig:BC},
but we will not discuss the latter. The partial unfolding is a
geometric object defined entirely by the combinatorial structure
of~$K$, and comes along with a canonical \mbox{projection~$p:\widehat{K}\to
K$.}

However, the partial unfolding~$\widehat{K}$ may not be a simplicial
complexes. In general~$\widehat{K}$ is a pseudo-simplicial complexes:
Let~$\Sigma$ be a collection of pairwise disjoint geometric simplices,
with simplicial attaching maps for some pairs
$(\s,\tau)\in\Sigma\times\Sigma$, mapping a subcomplex of~$\s$
isomorphically to a subcomplex of~$\tau$. Identifying the subcomplexes
accordingly yields the quotient space~$\Sigma/\mysim$, which is called
a \emph{pseudo-simplicial complex} if the quotient map
$\Sigma\to\Sigma/\mysim$ restricted to any $\s\in\Sigma$ is bijective.
The last condition ensures that there are no self-identifications
within each simplex~$\s\in\Sigma$.

\subsubsection*{The group of projectivities.}
Let~$\s$ and~$\tau$ be neighboring facets of a finite, pure simplicial
complex~$K$, that is, $\s\cap\tau$ is a ridge. Then there is exactly
one vertex in~$\s$ which is not a vertex of~$\tau$ and vice versa,
hence a natural bijection~$\langle\s,\tau\rangle$ between the vertex
sets of~$\s$ and~$\tau$ is given by
\[
  \langle\s,\tau\rangle:V(\s)\to V(\tau) :
  v\mapsto
  \begin{cases}
    v &\text{if} \quad v\in\s\cap\tau\\
    \tau\setminus\s &\text{if} \quad v=\s\setminus\tau.
  \end{cases}
\]
The bijection~$\langle\s,\tau\rangle$ is called a \emph{perspectivity}
from~$\s$ to~$\tau$.

\begin{figure}[t]
  \centering
  \begin{overpic}[width=.57\textwidth]{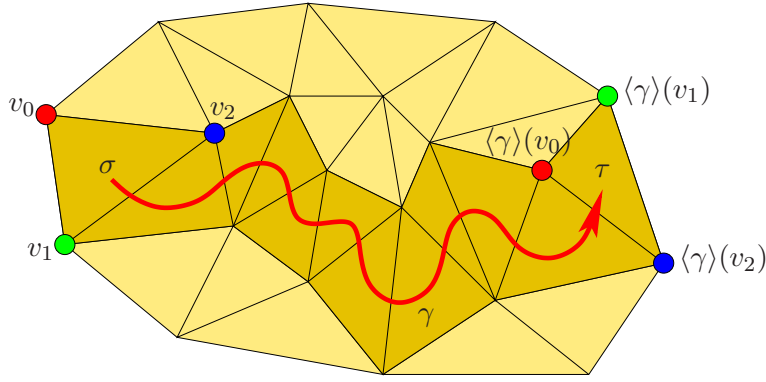}
    \put(13,32){$\mathbf{\s}$}
    \put(84.5,31.5){$\mathbf{\tau}$}
    \put(59,10.5){$\mathbf{\g}$}
    \put(0,41){$v_0$}
    \put(2.7,20){$v_1$}
    \put(29,40.5){$v_2$}
    \put(68.5,35.5){$\langle\g\rangle(v_0)$}
    \put(88.5,43){$\langle\g\rangle(v_1)$}
    \put(96.7,18.5){$\langle\g\rangle(v_2)$}
  \end{overpic}
  \caption{A projectivity from~$\s$ to~$\tau$ along the facet path~$\g$.\label{fig:projectivity}}
\end{figure}

A \emph{facet path} in~$K$ is a sequence $\g=(\s_0,\s_1,\dots,\s_k)$
of facets, such that the corresponding nodes in the dual graph~$\Gamma^*(K)$ form a
path, that is, $\s_i\cap\s_{i+1}$ is a ridge for all $0\leq i<k$; see
Figure~\ref{fig:projectivity}.  Now a
\emph{projectivity}~$\langle\g\rangle$ along~$\g$ is defined as the
composition of perspectivities $\langle\s_i,\s_{i+1}\rangle$,
thus~$\langle\g\rangle$ maps~$V(\s_0)$ to~$V(\s_k)$ bijectively via
\[
\langle\g\rangle = \langle\s_{k-1},\s_k\rangle\circ
\dots\circ\langle\s_1,\s_2\rangle\circ\langle\s_0,\s_1\rangle.
\]

We write $\g\,\g'=(\s_0,\s_1,\dots,\s_k,\dots,\s_{k+l})$ for
the \emph{concatenation} of two facet paths
$\g=(\s_0,\s_1,\dots,\s_k)$ and $\g'=(\s_k,\s_{k+1},\dots,\s_{k+l})$,
denote by $\g^-=(\s_k,\s_{k-1},\dots,\s_0)$ the \emph{inverse path}
of~$\g$, and we call~$\g$ a \emph{closed} facet path based at~$\s_0$
if $\s_0=\s_k$. The set of closed facet paths based at~$\s_0$ together
with the concatenation form a group, and a closed facet path~$\g$
based at~$\s_0$ acts on the set~$V(\s_0)$ via $\g\cdot
v=\langle\g\rangle(v)$ for $v\in V(\s_0)$. Via this action we obtain
the \emph{group of projectivities} $\Pi(K,\s_0)$ given by all
permutations $\langle\g\rangle$ of~$V(\s_0)$. The group of
projectivities is a subgroup of the symmetric group~$\Sym(V(\s_0))$ on
the vertices of~$\s_0$.

The projectivities along null-homotopic closed facet paths based
at~$\s_0$ generate the subgroup $\Pi_0(K,\s_0)\subgr \Pi(K,\s_0)$,
which is called the \emph{reduced group of projectivities}.  Finally,
if~$K$ is strongly connected then $\Pi(K,\s_0)$ and $\Pi(K,\s'_0)$,
respectively $\Pi_0(K,\s_0)$ and $\Pi_0(K,\s'_0)$, are isomorphic for
any two facets $\s_0,\s'_0\in K$. In this case we usually omit the
base facet in the notation of the (reduced) group of projectivities,
and write $\Pi(K)=\Pi(K,\s_0)$, respectively $\Pi_0(K)=\Pi_0(K,\s_0)$.

\subsubsection*{The odd subcomplex.}
Let~$K$ be locally strongly connected; in particular,~$K$ is pure. The
link of a co-dimension 2-face~$f$ is a graph which is connected
since~$K$ is locally strongly connected, and~$f$ is called \emph{even}
if the link $\Link_K(f)$ of~$f$ is bipartite, and \emph{odd} otherwise. We define the
\emph{odd subcomplex} of~$K$ as all odd co-dimension 2-faces (together
with their proper faces), and denote it by~$K_{\odd}$ (or sometimes~$\odd(K)$).

Assume that~$K$ is pure and admits a $(d+1)$-\emph{coloring} of its
graph~$\Gamma(K)$, that is, we
assign one color of a set of~$d+1$ colors to each vertex
of~$\Gamma(K)$ such that the two vertices of any edge carry different
colors. Observe that the $(d+1)$-coloring of~$K$ is minimal with
respect to the number of colors, and is unique up to renaming the
colors if~$K$ is strongly connected. Simplicial complexes that are
$(d+1)$-colorable are called \emph{foldable}, since a $(d+1)$-coloring
defines a non-degenerated simplicial map of~$K$ to the
$(d+1)$-simplex. In other places in the literature foldable simplicial
complexes are sometimes called balanced.

\begin{lem}\label{lem:foldable_and_odd}
  The odd subcomplex of a foldable
  simplicial complex~$K$ is empty, and
  the group of projectivities $\Pi(K,\s_0)$ is trivial. In particular we have
  $\langle\alpha\rangle=\langle\beta\rangle$ for any two facet paths~$\alpha$
  and~$\beta$ from~$\s$ to~$\tau$ for any two facets $\s,\tau\in K$.
\end{lem}



We leave the proof to the reader. As we will see in Theorem~\ref{thm:IJ_1} the odd subcomplex is of interest in particular for
its relation to $\Pi_0(K,\s_0)$ of a nice simplicial complex~$K$.
A projectivity \emph{around} an odd co-dimension 2-face~$f$ is a projectivity along a
facet path $\g\,l\,\g^-$, where~$l$ is a closed facet path in
$\Star_K(f)$ based at some facet $\s\in\Star_K(f)$, and~$\g$ is
a facet path from~$\s_0$ to~$\s$. The path $\g\,l\,\g^-$ is
null-homotopic since~$K$ is locally strongly simply
connected.

\begin{thm}[{Izmestiev \& Joswig~\cite[Theorem~3.2.2]{izmestiev_joswig:BC}}]
  \label{thm:IJ_1}
  The reduced group of projectivities~$\Pi_0(K,\s_0)$ of a
  nice
  simplicial complex~$K$ is generated by projectivities around the odd
  co-dimension 2-faces. In particular,~$\Pi_0(K,\s_0)$ is generated by
  transpositions.
\end{thm}

Consider a geometric
realization~$|K|$ of~$K$. To a given facet
path~$\g=(\s_0,\s_1,\dots,\s_k)$ in~$K$ we
associate a (piecewise linear) path~$\overline{\g}$ in~$|K|$ by
connecting the barycenter of~$\s_i$ to the barycenters of
$\s_i\cap\s_{i-1}$ and $\s_i\cap\s_{i+1}$ by a straight line for
$1\leq i<k$, and connecting the barycenters of~$\s_0$ and
$\s_0\cap\s_1$, respectively $\s_k$ and $\s_k\cap\s_{k-1}$. The
fundamental group~$\pi_1(|K|\setminus|K_{\odd}|,y_0)$ of a nice
simplicial complex~$K$ is generated by paths~$\overline{\g}$,
where~$\g$ is a closed facet path based at~$\s_0$, and~$y_0$ is the
barycenter of~$\s_0$;
see~\cite[Proposition~A.2.1]{izmestiev_joswig:BC}. Furthermore, due to
Theorem~\ref{thm:IJ_1} we have the group homomorphism
\begin{equation}\label{equ:h}
  \h_K:\pi_1(|K|\setminus|K_{\odd}|,y_0)\to\Pi(K,\s_0):[\overline{\g}]\mapsto\langle\g\rangle,
\end{equation}
where $[\overline{\g}]$ is the homotopy class of the
path~$\overline{\g}$ corresponding to a facet path~$\g$.

\subsubsection*{The partial unfolding}
\begin{figure}[htbp]
  \centering
  \begin{overpic}[width=.6\textwidth]{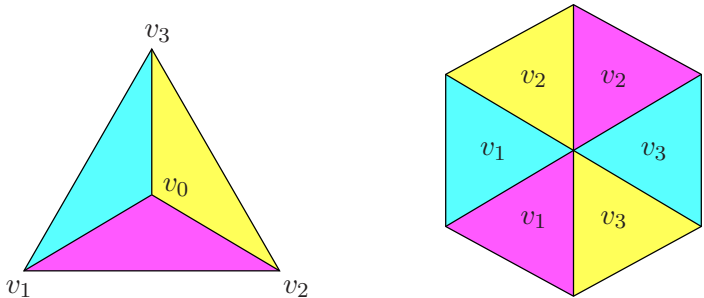}
    \put(22.5,18){$v_0$}
    \put(1,4){$v_1$}
    \put(39,4){$v_2$}
    \put(20,39.3){$v_3$}
    \put(66,23){$v_1$}
    \put(71.5,33){$v_2$}
    \put(88,23){$v_3$}
    \put(71.5,13.5){$v_1$}
    \put(82.5,33){$v_2$}
    \put(82.5,13.5){$v_3$}
  \end{overpic}
  \caption{The starred triangle and its partial unfolding. The complex on the
    right is the non-trivial connected component
    of the partial unfolding, indicated by the labeling of the facets
    by the vertices~$v_1$,~$v_2$, and~$v_3$. The second connected
    component is a copy of the starred triangle with all facets
    labeled~$v_0$; see also Example~\ref{exmp:branched_cover} for~$k=2$.
    \label{fig:part_unf}}
\end{figure}

Let~$K$ be a pure simplicial $d$-complex and set~$\Sigma$ as the set of all pairs $(|\s|,v)$,
where~$\s\in K$ is a facet and~$v\in\s$ is a vertex.
Thus each pair $(|\s|,v)\in\Sigma$ is a copy of the geometric simplex~$|\s|$
labeled by one of its vertices.
For neighboring facets~$\s$ and~$\tau$
of~$K$ we define the equivalence relation~$\sim$ by attaching
$(|\s|,v)\in\Sigma$ and $(|\tau|,w)\in\Sigma$ along their
common ridge~$|\s\cap\tau|$ if $\langle\s,\tau\rangle(v)=w$ holds. Now
the \emph{partial unfolding}~$\widehat{K}$
is defined as the quotient
space $\widehat{K} = \Sigma/\mysim$. The projection $p:\widehat{K}\to K$ is
given by the factorization of the map $\Sigma\to K:(|\s|,v)\mapsto
\s$; see Figure~\ref{fig:part_unf}.

The partial unfolding of a strongly connected simplicial complex is not
strongly connected in general. We denote by $\widehat{K}_{(|\s|,v)}$ the
connected component containing the labeled facet $(|\s|,v)$. Clearly,
$\widehat{K}_{(|\s|,v)}=\widehat{K}_{(|\tau|,w)}$ holds if and only if
there is a facet path~$\g$ from~$\s$ to~$\tau$ in~$K$ with
$\langle\g\rangle(v)=w$. It follows that the connected components
of~$\widehat{K}$ correspond to the orbits of the action
of~$\Pi(K,\s_0)$ on~$V(\s_0)$.  Note that each connected component of
the partial unfolding is strongly connected and locally strongly
connected~\cite[Satz~3.2.2]{witte:DIPL}. Therefore we do not
distinguish between connected and strongly connected components of the
partial unfolding.

The problem that the partial unfolding~$\widehat{K}$ may not be a
simplicial complex can be addressed in several ways. Izmestiev~\&
Joswig~\cite{izmestiev_joswig:BC} suggest barycentric subdivision
of~$\widehat{K}$, or anti-prismatic subdivision of~$K$.  A more
efficient solution (with respect to the size of the resulting
triangulations) is given in~\cite{witte:DIPL}.

\subsection{The partial unfolding as a branched cover.}

As preliminaries to this section we state two theorems by
Fox~\cite{fox:CSWS} and Izmestiev~\& Joswig~\cite{izmestiev_joswig:BC}. Together
they imply that under the ``usual connectivity assumptions'' the
partial unfolding of a
simplicial complexes is indeed a branched cover as suggested in the
heading of this section. For simplicial complexes the analog of these
topological connectivity requirements are
nice simplicial complexes.

\begin{thm}[{Izmestiev \& Joswig~\cite[Theorem 3.3.2]{izmestiev_joswig:BC}}]
  Let~$K$ be a nice simplicial complex. Then the restriction of
  $p:\widehat{K}\to K$ to the preimage of the complement of the odd
  subcomplex is a covering.
\end{thm}

\begin{thm}[{Fox~\cite[p.~251]{fox:CSWS}; Izmestiev~\&
    Joswig~\cite[Proposition 4.1.2]{izmestiev_joswig:BC}}]
  Let ~$J$ and~$K$ be nice simplicial complexes and let $f:J\to K$ be a
  simplicial map. Then the map~$f$ is a simplicial branched cover if and only if
  \[
  \Codim K_{\Sing}\geq2.
  \]
\end{thm}

Since the partial unfolding of a nice simplicial complex is nice
Corollary~\ref{cor:unf_are_covers} follows immediately.

\begin{cor}\label{cor:unf_are_covers}
  Let~$K$ be a nice simplicial complex. The projection
  $p:\widehat{K}\to K$ is a branched cover with the odd
  subcomplex~$K_{\odd}$ as its branching set.
\end{cor}

For the rest of this section let~$K$ be a nice
simplicial complex and let~$y_0$ be the barycenter
of~$|\s_0|$. The projection $p:\widehat{K}\to
K$ is a branched cover with $K_{\Sing}=K_{\odd}$ by Corollary~\ref{cor:unf_are_covers}, and Izmestiev~\&
Joswig~\cite{izmestiev_joswig:BC} proved that there is a bijection $\imath:p^{-1}(y_0)\to V(\s_0)$
that induces a group isomorphism $\imath_*:\Sym(p^{-1}(y_0))\to\Sym(V(\s_0))$
such that the following Diagram~\eqref{equ:partial_unf} commutes.
\begin{equation}\label{equ:partial_unf}
\xymatrix{ 
  \pi_1(|K|\setminus|K_{\odd}|,y_0) \ar[dr]^{\h_K} \ar[d]^{\m_p} & \\
  \M_p \ar[r]_{{\imath_*}\hspace{.15cm}\restr_{\M_p}} & \Pi(S,\s_0) \\
}
\end{equation}

In the case that the action of~$\Pi(K,\s_0)$
on~$V(\s_0)$ has only one non-trivial orbit we refer to the
unique non-trivial connected component of~$\widehat{K}$ corresponding to the non-trivial orbit as
\emph{the} partial unfolding. Otherwise
fix a set of generators of $\pi_1(|K|\setminus|K_{\odd}|,y_0)$
corresponding to closed (facet) paths around odd co-dimension 2-faces, and possibly further
generators of $\pi_1(|K|,y_0)$. Now each odd co-dimension 2-face
corresponds to exactly one non-trivial orbit of the 
$\Pi(K,\s_0)$-action, and~$K_{\odd}$ decomposes correspondingly. In this
spirit we can think of the empty set as the odd subcomplex
corresponding to a trivial orbit.

Consider a nice
simplicial complex~$K$, and a branched
cover $r:X\to Z$. Assume that there is a
homomorphism of pairs $\varphi:(Z,Z_{\Sing})\to(|K|,|K_{\odd}|)$, that is, $\varphi:Z\to|K|$
is a homomorphism with $\varphi(Z_{\Sing})=|K_{\odd}|$.
Then Theorem~\ref{thm:partial_unf} gives sufficient conditions for
$p:\widehat{K}\to K$ and $r:X\to Z$ to be equivalent branched covers.
It is the key tool in the proof of the main Theorem~\ref{thm:ch_4mf_main}
in Section~\ref{sec:constructing_4mfs}.

\begin{thm}\label{thm:partial_unf}
  Let~$K$ be a nice simplicial complex, and let $r:X\to Z$ be a simple branched
  cover. Further assume that there is a homomorphism of pairs
  $\varphi:(Z,Z_{\Sing})\to(|K|,|K_{\odd}|)$, and let $z_0\in Z$ be a
  point such that $y_0=\varphi(z_0)$ is the barycenter of~$|\s_0|$ for
  some facet $\s_0\in K$. The branched covers $p:\widehat{K}\to K$ and
  $r:X\to Z$ are equivalent
  if there is
  a bijection $\iota:r^{-1}(z_0)\to V(\s_0)$ 
  that induces a group isomorphism $\iota_*:\M_r\to\Pi(K,\s_0)$
  such that the diagram
  \begin{equation}\label{equ:partial_unf_and_bc}
    \xymatrix{
      \pi_1(Z\setminus Z_{\Sing},z_0) \ar[d]^{\m_r} \ar[r]^{\varphi_*} &
      \pi_1(|K|\setminus|K_{\odd}|,y_0) \ar[d]^{\h_K}\\
      \M_r
      \ar[r]^{\iota_*} & 
      \Pi(K,\s_0)
    }
  \end{equation}
  commutes. In particular, we have $\widehat{K}\cong X$.
\end{thm}
  
\begin{proof}
  Corollary~\ref{cor:unf_are_covers} ensures that $p:\widehat{K}\to K$ is indeed a branched cover, and
  commutativity of Diagram~\eqref{equ:partial_unf} and  Diagram~\eqref{equ:partial_unf_and_bc}
  proves commutativity of their composition:
  \[
  \xymatrix{
    \pi_1(Z\setminus Z_{\Sing},z_0) \ar[d]^{\m_r} \ar[r]^{\varphi_*} &
    \pi_1(|K|\setminus|K_{\odd}|,y_0) \ar[d]^{\m_p}\\
    \M_r
    \ar[r]^{\imath_*^{-1}\circ\iota_*} &
    \M_p
  }
  \]
  Theorem~\ref{thm:bc} completes the proof.
\end{proof}


\section{Color Equivalence of Simplicial Complexes}
\label{sec:color_equivalence}

\noindent
Consider two nice simplicial complexes~$K$ and~$K'$. The partial
unfoldings of two homeomorphic simplicial complexes need not to be
homeomorphic in general.  Here we present sufficient criteria for
$p:\widehat{K}\to K$ and $p':\widehat{K'}\to K'$ to be equivalent
branched covers.  Assume $K\cong K'$ and that the odd
subcomplexes~$K_{\odd}$ and~$K'_{\odd}$ are \emph{equivalent}, that
is, there is a homeomorphism of pairs $\varphi:(|K|,|K_{\odd}|)\to
(|K'|,|K'_{\odd}|)$. Let $\s_0\in K$ be a facet, and~$y_0$ the
barycenter of~$\s_0$, and assume that the image $y'_0=\varphi(y_0)$ is
the barycenter of~$|\s'_0|$ for some facet $\s'_0\in K'$. Now~$K$
and~$K'$ are \emph{color equivalent} if there is a bijection
$\psi:V(\s_0)\to V(\s'_0)$, such that
\begin{equation}\label{equ:color_equivalent}
  \psi_*\circ\h_K=\h_{K'}\circ\varphi_*
\end{equation}
holds,
where the maps $\varphi_*:\pi_1(|K|\setminus|K_{\odd}|,y_0)\to\pi_1(|K'|\setminus|K'_{\odd}|,y'_0)$
and $\psi_*: \Sym(V(\s_0))\to \Sym(V(\s'_0))$ are the group isomorphisms induced by~$\varphi$ and~$\psi$, respectively. 
Observe that this is indeed an equivalence relation. The name ``color
equivalent'' suggests that the pairs $(K,K_{\odd})$ and
$(K',K'_{\odd})$ are equivalent, and that the ``colorings'' of~$K_{\odd}$ and~$K'_{\odd}$
by the $\Pi(K)$-action, respectively $\Pi(K')$-action, of projectivities
around odd faces are
equivalent. Proposition~\ref{prop:color_equiv} justifies this name.

\begin{prop}\label{prop:color_equiv}
  Let~$K$ and~$K'$ be color equivalent simplicial complexes. Then the
  branched covers $p:\widehat{K}\to K$ and $p':\widehat{K'}\to K'$ are
  equivalent.
\end{prop}

\begin{proof}
  With the notation of Equation~\eqref{equ:color_equivalent} we have that
  \[
  \xymatrix{
    & \pi_1(|K|\setminus|K_{\odd}|,y_0) \ar[r]^{\varphi_*}
    \ar[dl]_{\m_p} \ar[d]^{\h_K}&
    \pi_1(|K'|\setminus|K'_{\odd}|,y'_0)
    \ar[dr]^{\m_{p'}} \ar[d]_{\h_{K'}}& \\
    \M_p \ar[r]^{\imath_*} & 
    \Pi(K,\s_0)\ar[r]^{\psi_*} &
    \Pi(K',\s'_0)& 
    \M_{p'} \ar[l]_{\imath'_*}
  }
  \]
  commutes, since the Diagram~\eqref{equ:partial_unf} commutes and
  Equation~\eqref{equ:color_equivalent} holds. Theorem~\ref{thm:bc} completes the proof.
\end{proof}

\subsubsection*{The anti-prismatic subdivision.}
\begin{figure}[htbp]
  \centering
  \begin{overpic}[width=.8\textwidth]{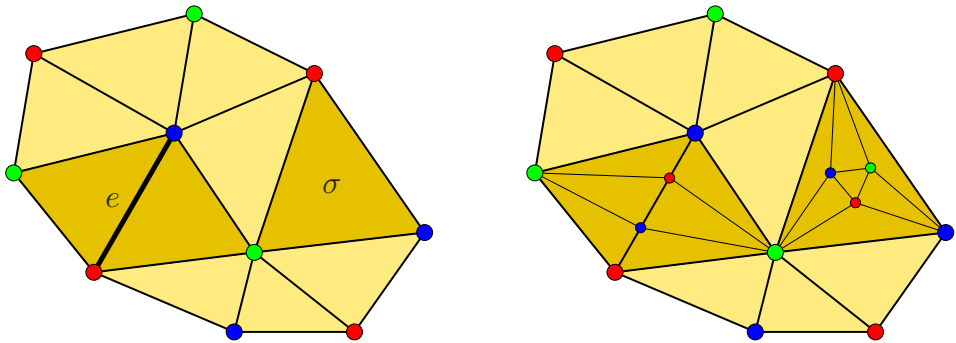}
    \put(11.5,15){\large $e$}
    \put(33.8,16.5){\large $\s$}
  \end{overpic}
  \caption{Anti-prismatic
    subdivision of the edge $e$ and the triangle~$\s$ of a foldable
    simplicial complex.
    \label{fig:anti_prism_f2}}
\end{figure}

Let~$c_k$ be the simplicial complex obtained from the boundary complex
of the $(k+1)$-dimensional cross polytope by removing one facet.
Alternatively, define~$c_k$ as the simplicial complex arising from the
Schlegel diagram of the $(k+1)$-dimensional cross polytope. To be more explicit, let
$\s=\{+v_i\}_{0\leq i\leq k}$ be the vertices of the $k$-simplex. Then
the facets of~$c_k$ are defined as all subsets $\s'\not=\s$ of $\{\pm
v_i\}_{0\leq i\leq k}$ such that either~$+v_i\in\s'$
\mbox{or~$-v_i\in\s'$} holds. The complex~$c_k$ and the $k$-simplex
are PL-homeomorphic with isomorphic boundaries, and~$c_k$ is
$(k+1)$-colorable by assigning the same color to~$+v_i$ and~$-v_i$, since
$\{+v_i,-v_i\}$ is not an edge. The \emph{anti-prismatic
  subdivision}~$a_f(K)$ of a $k$-face~$f$ of a simplicial
$d$-complex~$K$ is obtained from~$K$ by replacing~$\Star_K(f)$ by the
join of $c_k$ with~$\Link_K(f)$, that is
\[
a_f(K)=(K\setminus \Star_K(f))\cup(c_k*\Link_K(f)).
\]
See Figure~\ref{fig:anti_prism_f2} for a an example of an
anti-prismatic subdivision of an edge and a triangle of a foldable
simplicial complex.

The \emph{anti-prismatic subdivision}~$a(K)$ of a simplicial
complex~$K$ is defined by recursively anti-prismatically subdividing
all faces of~$ K$ from the facets down to the edges.  Observe
that~$a_f(K)$, and hence~$a(K)$, are PL-homeomorphic to~$K$, and
that~$a_f(K)$ and~$a(K)$ inherit niceness from~$K$.
For sake of brevity we omit the (lengthy but straight forward)
proof of
Proposition~\ref{prop:anti-prismatic}. The result is similar to
Proposition~A.1.3 and following in~\cite{izmestiev_joswig:BC}, and
explicit proofs can be found in~\cite[Section~1.3.1]{witte:DISS}.

\begin{prop}\label{prop:anti-prismatic}
  Let~$K$ be a nice simplicial complex and $f\in K$ a face. The simplicial complexes~$a_f(K)$,~$a(K)$ and~$K$ are color equivalent.
\end{prop}



\subsection{Prescribing the odd subcomplex.}

Theorem~\ref{thm:partial_unf} made it clear, that it is essential to
control the odd subcomplex if one tries to determine the behavior of
the partial unfolding as a branched cover,
e.g. in the construction of combinatorial 4-manifolds in Section~\ref{sec:constructing_4mfs}.

Let~$K$ be a strongly connected and foldable simplicial complex of
dimension~$d$, and fix a $(d+1)$-coloring using the colors
$[d+1]=\{0,1,\dots,d\}$. Then the
$\{i_0,i_1,\dots,i_k\}$-\emph{skeleton} is the subcomplex of~$K$
induced by the vertices colored $\{i_0,i_1,\dots,i_k\}$. Observe that
the $\{i_0,i_1,\dots,i_k\}$-skeleton is a pure simplicial complex of
dimension~$k$.

\begin{prop} \label{prop:odd_subcomplex}
  Let~$K$ be a foldable
  combinatorial manifold of dimension~$d$ and let~$F$ be a
  co-dimension 1-manifold (possibly with more than one connected
  component) embedded in the
  $\{i_0,i_1,\dots,i_{d-1}\}$-skeleton
  of~$K$.  Further assume that all facets (and their proper faces)
  of the boundary~$\partial F$ of~$F$ not contained entirely in $\partial K$, for short
  the \emph{closure}
  $\closure(\partial F\setminus \partial K)$,
  are embedded in the $\{i_0,i_1,\dots,i_{d-2}\}$-skeleton.
  Then~$\closure(\partial F \setminus \partial K)$ can
  be realized as the odd subcomplex of some simplicial complex~$K'$,
  that arises from~$K$ by stellar subdivision of edges in the
  $\{i_{d-1},i_d\}$-skeleton. The
  complex~$K'$ is $(d+2)$-colorable by extending the coloring of~$K$, and the odd subcomplex lies in the
  $\{i_0,i_1,\dots,i_{d-2}\}$-skeleton.
\end{prop}

\begin{proof}
  Every $(d-1)$-simplex
  in~$F$ has exactly one $i_{d-1}$-colored vertex since~$F$ is
  foldable. Hence the vertex stars of all $i_{d-1}$-colored
  vertices cover~$F$, that is,
  \begin{equation}\label{equ:vertex_stars_of_F}
    F = \bigcup_{\text{$v$ is $i_{d-1}$-colored}} \Star_F(v),
  \end{equation}
  and the vertex stars intersect in the
  $\{i_0,i_1,\dots,i_{d-2}\}$-skeleton. Further, a $(d-2)$-face~$g\in
  F$ (a ridge in~$F$) is contained in an odd number of vertex stars of
  $i_{d-1}$-colored vertices of~$F$ if and only if~$g \in \partial F$
  since~$F$ is an embedded combinatorial manifold.
  Observe that stellar subdivision of an edge~$e$ changes the parity of
  $\Link_K(g)$ of each co-dimension 2-face $g\in\Link_K(e)\setminus
  \partial K$. The odd subcomplex resulting from a series of stellar
  subdivisions of edges is the symmetric difference of the edge links.

  Since~$K$ is a combinatorial manifold the vertex star~$\Star_K(v)$
  of an $i_{d-1}$-colored vertex~$v\in F$ is a $d$-ball, which is the
  join of~$v$ with an $(i_0,i_1,\dots,i_{d-2},i_d)$-colored
  $(d-1)$-ball if~$v\in\partial K$, and which is the join with an
  $(i_0,i_1,\dots,i_{d-2},i_d)$-colored $(d-1)$-sphere otherwise. The
  vertex star~$\Star_F(v)$ divides~$\Star_K(v)$ into two connected
  components, and we will call these two connected components of
  $|\Star_K(v)|\setminus|\Star_F(v)|$ the two \emph{sides}
  of~$\Star_F(v)$, mimicking the topological concept of a two-sided
  manifold (embedded in an orientable space); see
  Figure~\ref{fig:star_of_v} for a 3-dimensional example. The
  link~$\Link_K(\{v,w\})$ of an $\{i_{d-1},i_d\}$-colored
  edge~$\{v,w\} \in \Star_K(v)$ is a $(d-2)$-sphere in the
  $\{i_0,i_1,\dots,i_{d-2}\}$-skeleton of~$\partial\Star_K(v)$.
  Moreover, the vertex stars of all $\{i_{d-1},i_d\}$-colored edges
  $\{v,w\}\in\Star_K(v)$ cover~$\Link_K(v)$.  Thus if we stellar
  subdivide all $\{i_{d-1},i_d\}$-edges in one side of~$\Star_F(v)$ we
  obtain~$\Link_F(v)$ as the odd subcomplex.

  \begin{figure}[t]
    \center
    \psfrag{0}{\small\;1}
    \psfrag{1}{\small\;0}
    \psfrag{2}{\small\;2}
    \psfrag{3}{\small\;3}
    \begin{overpic}[width=.38\textwidth]{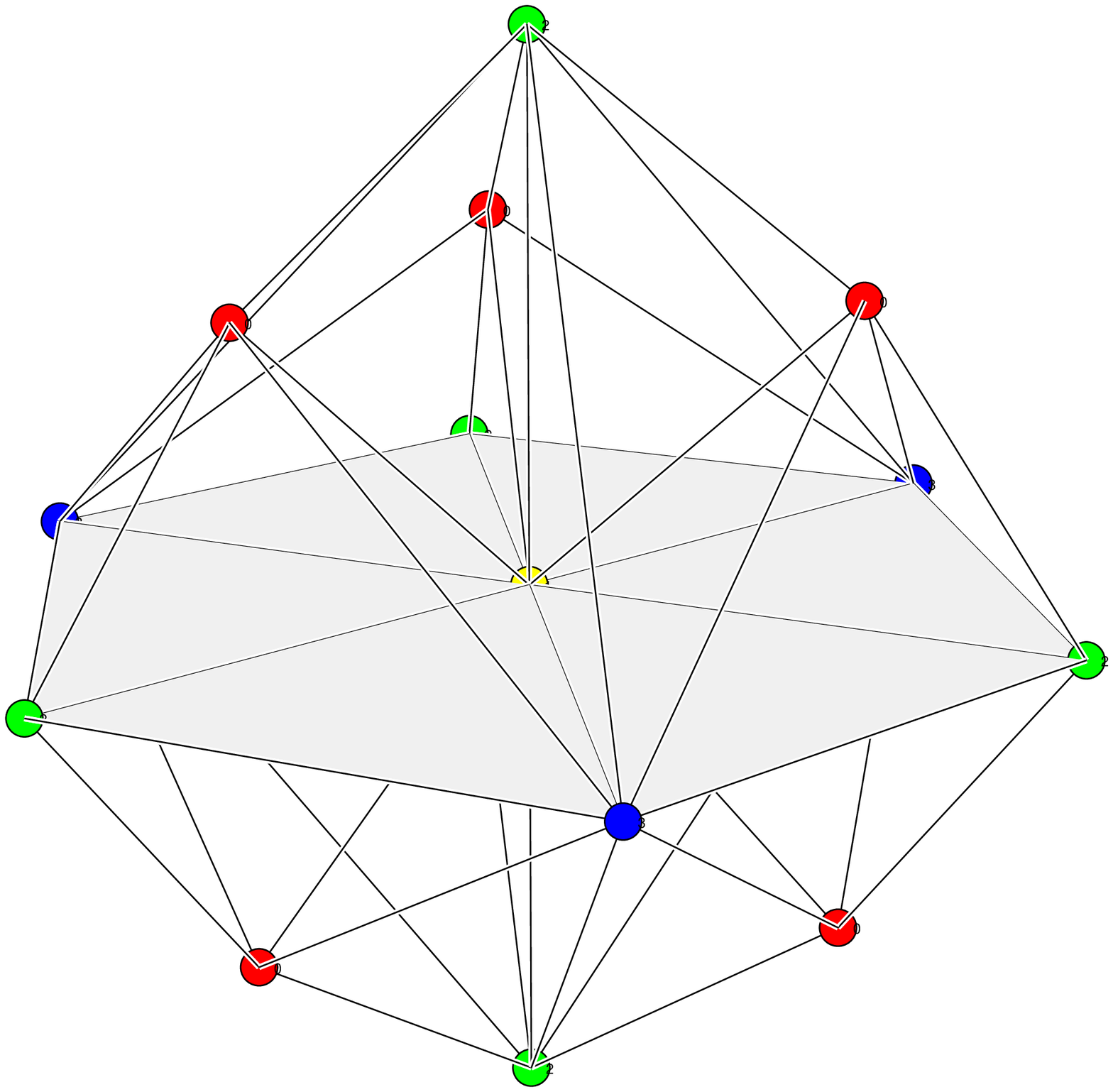} 
      \put(34.5,33){$F$}
    \end{overpic}
    \hspace{.9cm}
    \begin{overpic}[width=.38\textwidth]{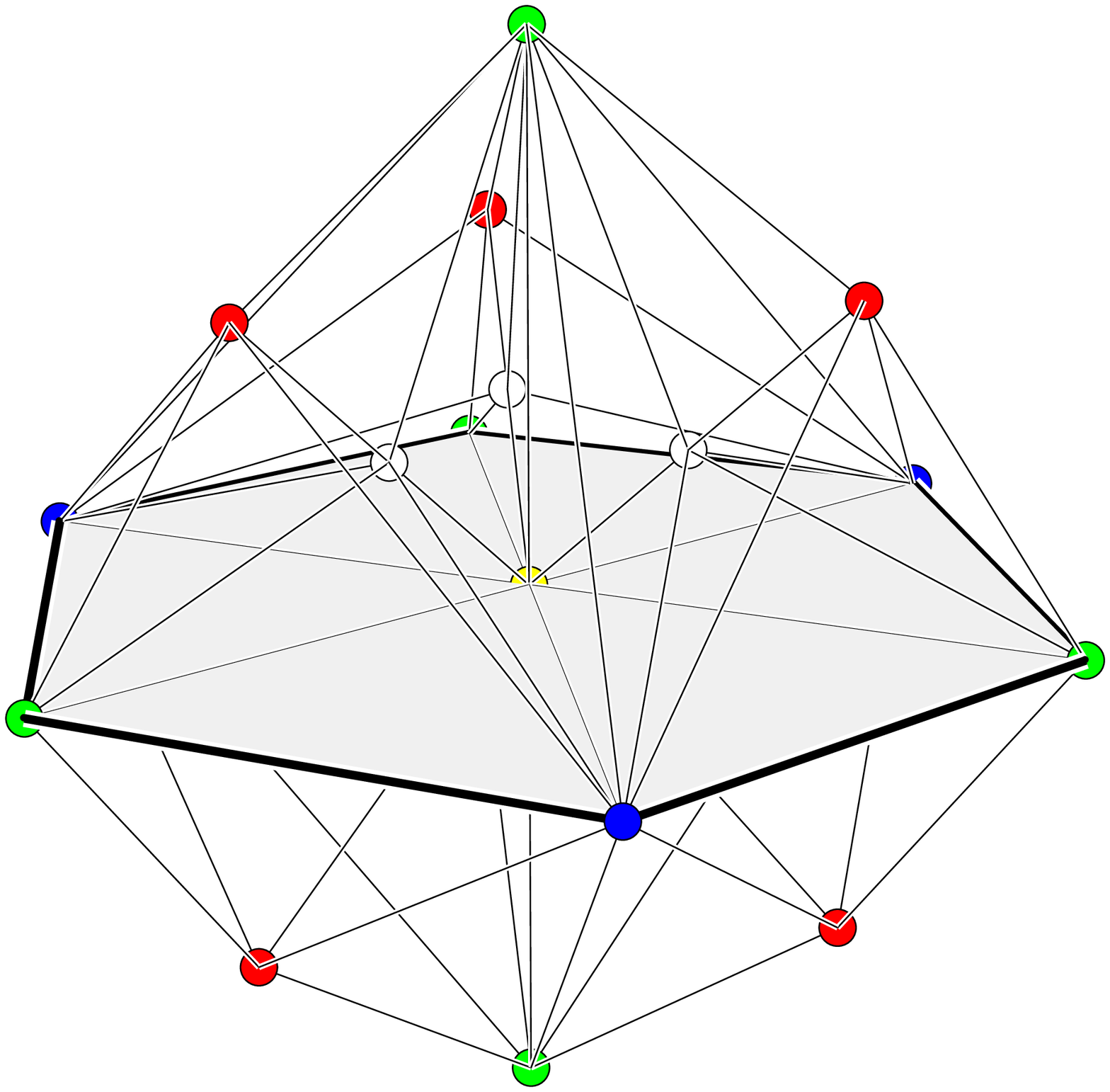} 
      \put(34.5,33){$F$}
    \end{overpic}
    \caption{Vertex star of an 0-colored vertex~$v\in
      F$, on the right after stellar subdivisions of all
      $\{0,1\}$-edges in one side of $\Star_F(v)$. The parity of
      the edges in $\Link_F(v)$ changes.\label{fig:star_of_v}}
  \end{figure}

  Finally we construct the desired odd subcomplex $\closure(\partial F
  \setminus \partial K)$ as the symmetric difference of vertex links of
  all $i_{d-1}$-colored vertices in~$F$.

  The resulting complex~$K'$ is $(d+2)$-colorable by assigning a new
  color to the vertices introduced by stellar subdivision of
  edges. If an edge~$e$ is subdivided twice, use the original colors
  of~$e$ to color the two new vertices.
\end{proof}

\begin{rem}\label{rem:odd_subcomplex}
  Observe that a projectivity based at~$\s_0$ around an odd co-dimension
  2-face constructed via Proposition~\ref{prop:odd_subcomplex} exchanges
  the two vertices of~$\s_0$ colored~$i_{d-1}$ and~$i_d$.
\end{rem}  

We conclude this section with a characterization of some co-dimension 2-manifolds
which by Proposition~\ref{prop:odd_subcomplex} can be
realized as an odd subcomplex in the
$\{i_0,i_1,\dots,i_{d-2}\}$-skeleton.

\begin{lem}\label{lem:odd_subcomplex}
  Let~$|K|$ be some geometric realization of a foldable combinatorial
  $d$-manifold~$K$. An orientable PL $(d-1)$-manifold~$F$ may be embedded
  in the $\{i_0,i_1,\dots,i_{d-2},d-1\}$-skeleton of (a refinement
  of)~$K$ with boundary~$\partial F$ embedded in the
  $\{i_0,i_1,\dots,i_{d-2}\}$-skeleton if there is an embedding
  $F\to|K|$. Note that the last color in the coloring of
  the embedding of~$F$ is explicitly required to be~$d-1$.
\end{lem}

\begin{proof}
  Simplicial approximation of the embedding $F\to|K|$ yields an embedding of~$F$
  in the co-dimension 1-skeleton of some refinement~$K'$ of~$K$. Let~$b(K')$
  be the barycentric subdivision
  of~$K'$ with each vertex colored by
  the dimension of its originating face. The embedding $F\to K'$
  yields an embedding~$\imath: F\to b(K')$ of~$F$ in the $\{0,1,\dots,d-1\}$-skeleton
  of~$b(K')$, with~$\partial F$ embedded in the
  $\{0,1,\dots,d-2\}$-skeleton. Further we have that the vertex stars of
  all $(d-1)$-colored vertices cover~$\imath(F)$; see
  Equation~\eqref{equ:vertex_stars_of_F}.
  
  It remains to show, how to
  ``push''~$F$ into the desired skeleton. The
  $\{0,1,\dots,d-2,d-1\}$-skeleton of~$b(K')$ differs from the
  $\{i_0,i_1,\dots,i_{d-2},d-1\}$-skeleton by one
  color~$c=\{0,1,\dots,d-2\}\setminus\{i_0,i_1,\dots,i_{d-2}\}$, that
  is, replacing~$c$ by~$d$ in $\{0,1,\dots,d-2,d-1\}$ yields $\{i_0,i_1,\dots,i_{d-2},d-1\}$.  
  For each $(d-1)$-colored vertex
  $v\in \imath(F)$ choose one of the two sides of $\Star_{\imath(F)}(v)$
  consistent with the orientation of~$F$.
  This may
  be done since~$F$ is orientable. Let~$v\in \imath(F)$ be $(d-1)$-colored,
  let~$D_v$ be the chosen side of $\Star_{\imath(F)}(v)$,
  and let~$V_c$ denote the set of all
  $c$-colored vertices in $\Link_{\imath(F)}(v)$; see Figure~\ref{fig:star_of_v}. Now we obtain the desired
  embedding $\imath':F\to b(K')$ by replacing $\Star_{\imath(F)}(v)$ with
  \[
  \bigcup_{w\in V_c} v * (\Link_{b(K')}(\{v,w\})\cap D_v)\cong \B^{d-1}.
  \]
  Here it is important that the triangulation of~$b(K')$ may have to
  be refined further by anti-prismatic subdivision.  The map
  $\imath':F\to b(K')$ is an embedding of~$F$ since we replace
  $(d-1)$-balls by $(d-1)$-balls, and two $(d-1)$-balls
  in~$\imath'(F)$ intersect as in~$\imath(F)$ due to the consistent
  choice of the sides of $\Star_{\imath(F)}(v)$.
\end{proof}

\subsection{Extending triangulations.}

We present a technique how to extend a partial triangulation of some
topological space to the entire space (e.g. a triangulation
of~$\Sph^{d-1}$ to~$\B^d$) while considering certain restraints on the
colorability of the triangulations. This technique is crucial in the
constructions in Section~\ref{sec:constructing_4mfs}.

A first assault on how to extend triangulation and coloring is by Goodman~\&
Onishi~\cite{goodman_onishi:ETC}, who proved that a 4-colorable
triangulation of~$\Sph^2$ may be extended to a 4-colorable
triangulation of~$\B^3$. Their result was improved independently by
Izmestiev~\cite{izmestiev:EOC} and~\cite{witte:DISS,witte:DIPL} to arbitrary
dimensions. In~\cite[Theorem~2.3]{witte:DISS} further restrictions (e.g. regularity) are
taken into account.

\begin{lem} \label{lem:extend_coloring}
  Let~$S$ be a $k$-colored combinatorial $(d-1)$-sphere. Then there
  exists a combinatorial $d$-ball~$B$ with boundary~$\partial B$ equal
  to~$S$ such that the $k$-coloring of~$S$ may be extended to a
  $\max\{k,d+1\}$-coloring of~$B$.
\end{lem}

\begin{figure}[htbp]
  \centering
  \includegraphics[width=.76\textwidth]{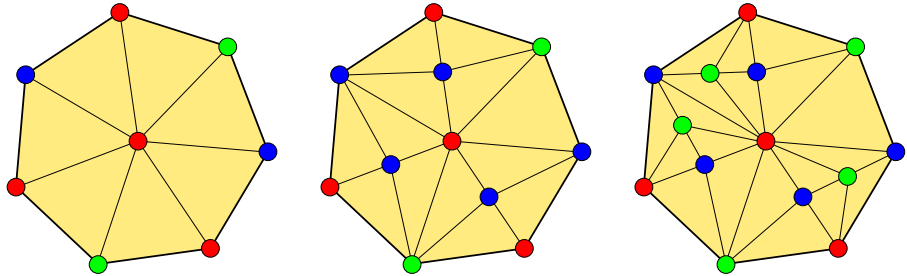}
  \caption{
    Extending a triangulation
    of~$\Sph^1$ and its 3-coloring to~$\B^2$.
    \label{fig:ext_coloring0}}
\end{figure}

We sketch the construction in~\cite{witte:DISS}.
Set $B_0 = a*S$ as the cone
over~$S$ with apex~$a$. We have
$k\geq d$ and in the case $k=d$ set $B = B_0$ and color~$a$ with a new
color. Otherwise fix an ordering $c_0<c_1<\dots<c_d$ of the first~$d+1$
colors in the coloring of~$S$ and color~$a$ by~$c_0$. For $1\leq i\leq d$ we obtain~$B_i$ from~$B_{i-1}$ by
stellarly subdividing all edges with both vertices colored~$c_{i-1}$. The
new vertices are colored~$c_i$. Theorem~2.3 in~\cite{witte:DISS}
ensures that $B=B_d$ is properly colored; see Figure~\ref{fig:ext_coloring0}.


Let~$X=X^d$ be a finite CW-complex of dimension~$d$ with $l$-cells
$\{e^l_\alpha\}_\alpha$, closed cells
$\{C^l_\alpha\}_\alpha=\{\closure(e^l_\alpha)\}_\alpha$, attaching
maps $\varphi_\alpha:\partial C^l_\alpha \to X^{l-1}$, and $l$-skeleton
\[
X^l=X^{l-1} \dotcup \big(\bigcup_\alpha \hspace{-.33cm}\cdot\hspace{.24cm} e^l_\alpha\big).
\]
A finite CW-complex is \emph{regular} if the attaching maps
$\varphi_\alpha:\partial C^l_\alpha\to X^{l-1}$ (restricted to their image) are homeomorphisms for
all $1\leq l\leq d$; see Hatcher~\cite[p.~5]{hatcher:AT}. We call a
simplicial complex $K\cong X$ a \emph{triangulation} of~$X$, if~$K$
refines the cell structure of~$X$, that is, the $(d-1)$-skeleton
of~$K$ is a triangulation of the CW-complex~$X^{d-1}$.

A subset $Y\subset \{e^l_\alpha\}_\alpha$ is called a
\emph{subcomplex} if for each closed cell $C^l_\alpha\in Y$ all cells
in the image of $\varphi_\alpha:\partial C^l_\alpha\to X^{l-1}$ are also
in~$Y$. Hence~$Y$ is also a CW-complex, and~$Y$ is regular if~$X$ is
regular.  For example, any $l$-skeleton~$X^l$ is a subcomplex of~$X$.
We call a triangulation of a subcomplex~$Y\subset X$ a \emph{partial
  triangulation} of~$X$.

\begin{prop}\label{prop:extend_cw_complex}
  Let~$X$ be a finite regular CW-complex of dimension at most~$4$, and
  let~$Y\subset X$ be a subcomplex. Then any triangulation and $k$-coloring
  of~$Y^l$ can be extended to a triangulation and $\max\{k,l+1\}$-coloring of~$X^l$.
\end{prop}

\begin{proof}
  We prove by induction on $0\leq i \leq l$ that there exists a
  triangulation of the $i$-skeleton~$X^i$ which can be colored
  with $\max\{k,i+1\}$ colors such that the triangulation and coloring
  of~$X^i$ extend the triangulation and coloring of~$Y^i$. This clearly
  holds for~$i=0$, and for~$i=l$ we obtain Proposition~\ref{prop:extend_cw_complex}.

  Let~$i\geq1$ and let~$e^i_\alpha$ be an $i$-cell of~$X^i$ not
  contained in~$Y^i$. By
  induction~$X^{i-1}$ is triangulated and colored using
  $\max\{k,i+1\}$ colors, and the triangulation of~$X^{i-1}$
  extends triangulation and coloring of~$Y^{i-1}$.
  Since~$X$ is regular, the image of the attaching map
  $\varphi_\alpha:\partial C^i_\alpha\to X^{i-1}$ is a $(i-1)$-sphere
  with a triangulation induced by the triangulation of~$X^{i-1}$. Since $i\leq d$ is at
  most~4, every simplicial $(i-1)$-sphere is a combinatorial $(i-1)$-sphere.
  Now Lemma~\ref{lem:extend_coloring} extends triangulation and
  coloring of $\varphi_\alpha(\partial C^i_\alpha)$ to the entire $i$-ball~$C^i_\alpha$. Since the $i$-balls
  $\{C^i_\alpha\}_\alpha$ intersect pairwise only in~$X^{i-1}$, extending the triangulation of~$\partial C^i_\alpha$ to its interior for each
  $i$-cell~$e^i_\alpha$ yields the desired triangulation of~$X^i$.
\end{proof}

\begin{rem} A similar result to
  Propositions~\ref{prop:extend_cw_complex} for a partial
  triangulation~$R$ of a relative handlebody decomposition
  $|R|=N^{-1}\subset N^0 \subset \dots \subset N^4 = N$
  of the pair $(|R|,N)$ can be found in~\cite[Proposition~2.12]{witte:DISS}.
\end{rem}

\begin{rem}
  Proposition~\ref{prop:extend_cw_complex} is not applicable in higher
  dimensions. For example, let~$H$ be a triangulation of the
  Poincar\'e homology sphere; see Bj\"orner~\&
  Lutz~\cite{bjoerner_lutz:PHS,bjoerner_lutz:EGM-PHS}
  and~\cite{witte:DIPL}. The double suspension $\Susp^2(H)$ is
  homeomorphic to~$\Sph^5$~\cite{MR554220,cannon:SCLDM}, yet not a combinatorial sphere: There are
  two vertices with $\Susp(H)\not\cong\Sph^4$ as vertex links. Consider
  the cell decomposition of the 6-ball given by the triangulation
  $\Susp^2(H)$ of~$\Sph^5$ plus an additional 6-cell. Now,
  when attaching the final
  6-cell, one can not apply the inductive argument for the two
  vertices with $\Susp(H)$ as vertex links.
\end{rem}


\section{Constructing Combinatorial 4-Manifolds}
\label{sec:constructing_4mfs}

\noindent
The main result, the construction of a simplicial
4-sphere~$S$ such that the partial unfolding~$\widehat{S}$ is
PL-homeomorphic to a given closed oriented PL 4-manifold~$M$, is developed in
Section~\ref{sec:4mf_as_unfoldings}.
Prior to giving a combinatorial construction of~$M$,
we will review the topological situation.

Let~$W^3$ be a 3-manifold. Following
Montesinos~\cite{montesinos:MIC}, we call two given branched
covers
$p_1, p_2: W^3 \to \mathbb{S}^3$ branched over
links~$L_1$ and~$L_2$, respectively, \emph{cobordant} if
there exists a branched cover
$p: W^3 \times [0,1] \to \mathbb{S}^3 \times [0,1]$
which is equal to~$p_1$ if restricted to $W^3 \times \{0\}$, and equal to~$p_1$ if restricted to
$W^3 \times \{1\}$, and is branched over an immersed PL 2-manifold
with a
boundary equal to the disjoint union $L_1 \dotcup L_2$. The branched
cover~$p$ is called a \emph{cobordism}.

A (surprisingly) useful technique is to attach a \emph{trivial sheet}.
Given a $k$-fold branched cover $p:X\to\B^4$, (with sheets numbered
$0,1,\dots,k-1$) we want to add another sheet to the covering without
changing the topology of the covering space~$X$. To this end add a
2-ball~$D$ to the branching set of~$p$ such that~$\partial D$ is
contained in~$\partial\B^4$, and let a meridial loop around~$D$
correspond to the transposition $(1,k)$. The covering space~$X'$ of
the branched cover obtained this way is the union of~$X$ and a 4-ball
attached to~$\partial X$ along a 3-ball, thus~$X'\cong X$ holds.



\subsection{4-Manifolds as branched covers}
\label{sec:4mf_as_branched_covers}

Every closed oriented PL 4-manifold admits a handle representation
\[
M=H^0 \cup \lambda H^1 \cup \mu H^2 \cup \g H^3 \cup H^4;
\]
see~\cite[Example~2.2]{witte:DISS}. With $M_A = H^0 \cup \lambda H^1
\cup \mu H^2$, and $M_B = H^0 \cup \g H^1$ by duality, we obtain~$M$
as the union $M_A \cup_h M_B$, where~$h$ is the attaching map. That
is, we paste~$M_A$ and~$M_B$ together along their common boundary $\g
\;\sharp\; (\Sph^1 \times \Sph^2$), the connected sum of~$\g$ copies of
$\Sph^1 \times \Sph^2$.  In fact, Montesinos~\cite{montesinos:HDCM}
proved that $H^0 \cup \lambda H^1 \cup \mu H^2$ alone topologically
determines~$M$.

The construction of the branched cover $r: M \to \Sph^4$
proceeds in two steps.
In the first step (see Montesinos~\cite{montesinos:CSR})
the 4-manifolds~$M_A$ and~$M_B$ are constructed as simple 3-fold branched
covers~$r_A$ and~$r_B$ of the 4-ball~$\B^4$. Since $M_B = H^0
\cup \g H^1$ is of the form $H^0 \cup \lambda H^1 \cup \mu H^2$ it
suffices to show how to construct~$r_A:M_A\to\B^4$.

Although $\partial M_A = \partial M_B$ holds, the branching sets of~$r_A$ and~$r_B$ restricted to the common
boundary $\g \; \sharp \; (\Sph^1 \times \Sph^2)$ of~$M_A$ and~$M_B$ may not be equivalent, and
$M_A \cup_h M_B$ is not the covering space of a branched cover
$r: M \to \Sph^4$ in general. Hence in the second step we construct
a cobordism between~$r_A\left|_{\partial M_A}\right.$
and~$r_B\left|_{\partial M_B}\right.$, that is, a branched
cover $r_H : H \to \Sph^3 \times [0,1]$ with covering space $H \cong (
\g \; \sharp \; (\Sph^1 \times \Sph^2) ) \times [0,1]$ which satisfies
\begin{equation}\label{equ:cobordism}
r_H \restr_{\left( \g \; \sharp \; (\Sph^1 \times \Sph^2) \right) \times \{0\}}
= r_A \restr_{\partial M_A} \quad \text{and} \quad
r_H \restr_{\left( \g \; \sharp \; (\Sph^1 \times \Sph^2) \right) \times \{1\}}
= r_B \restr_{\partial M_B}.
\end{equation}
The cobordism~$r_H$ is branched over a PL 2-manifold with
a finite number of singularities and boundary equal to the disjoint union of the branching sets of
$r_A\left|_{\partial M_A}\right.$ and $r_B\left|_{\partial
    M_B}\right.$. The boundary of the
covering space~$H$ is homeomorphic to two disjoint copies of $\g \;
\sharp \; (\Sph^1 \times \Sph^2)$, and we have
$M \cong M_A \cup_{\Id} H \cup_{\Id} M_B$. Note that~$r_H$ is a 4-fold
cover in general, thus we must add a fourth, trivial sheet
to~$r_A$, respectively~$r_B$.

The existence of such a cobordism, and hence the representation of~$M$
as a branched cover of~$\Sph^4$, was first observed by
Piergallini~\cite{piergallini:FMS4}. 
The following diagram illustrates this approach.
\[
\xymatrix{
  M \ar[d]^{r} & \cong & M_A \ar[d]^{r_A} & \cup_{\Id} & H \ar[d]^{r_H} & \cup_{\Id} & M_B
  \ar[d]^{r_B}\\
  \Sph^4 & \cong & \B^4 & \cup_{\Id} & \Sph^3 \times [0,1] & \cup_{\Id} & \B^4
}
\]


\subsubsection*{Construction of $M_A$.}
In the following we will sketch a construction of $r_A: M_A \to \B^4$ as a 3-fold
branched cover branched over a ribbon manifold. This
construction is due to Montesinos, and we omit the proofs; refer
to~\cite{montesinos:CSR} for further details. 

First, consider the 4-manifold $W = H^0 \cup \lambda H^1$
which consists of a single 4-ball and 1-handles
only. It arises as the \emph{standard branched cover}
$r_W: W \to \B^4$ branched along $\lambda+2$ unlinked and unknotted copies
of~$\B^2$.
Let~$u:\RR^4\to\RR^4$ be the reflection on the hyperplane given by $x_1=0$, that
is,~$u$ maps $(x_1,x_2,x_3,x_4)$ to $(-x_1,x_2,x_3,x_4)$.  The
covering space~$W$ is obtained from
$[-1,1]^3\times[-1,2]$ by the following identifications on
its boundary. Consider the subset~$\mathcal{A}$ of
$[-1,1]^3$ consisting of~$2\lambda$ disjoint rectangles
given by
\[
\mathcal{A} = \bigcup_{i=1}^\lambda
\SetOf{(x_1,x_2,x_3)\in[-1,1]^3}{
  x_3=1 \;\text{and}\;
  x_1 \in \pm \left[\frac{2i-1}{2\lambda+1},\frac{2i}{2\lambda+1}\right]}.
\]

Now identify a point~$x \in [-1,1]^3\times[-1,2]$ with its
image~$u(x)$ if~$x$ lies in the top or bottom facet, that is, if
$x_4=-1$ or $x_4=2$, or if $(x_1,x_2,x_3) \in \mathcal{A}$ and
$x_4\in[-\frac{1}{2},\frac{1}{2}]\cup[\frac{3}{2},2]$.
If~$\mathcal{R}$ denotes the equivalence relation given by these
identifications, we have $([-1,1]^3\times[-1,2])/\mathcal{R}\cong W$.

\begin{figure}[t]
  \center
  \psfrag{[0,1]xB3}{$([-1,1]^2\times[0,1])/\mathcal{R}'$}
  \psfrag{[-1,2]xB3}{$([-1,1]^2\times[-1,2])/\mathcal{R}$}
  \psfrag{u}{$u$}
  \psfrag{pW}{$r_W$}
  \psfrag{x}{$x_1$}
  \psfrag{y}{$x_3$}
  \psfrag{z}{$x_4$}
  \psfrag{-1}{-1}
  \psfrag{0}{0}
  \psfrag{1}{1}
  \psfrag{2}{2}
  \includegraphics[width=.7\linewidth]{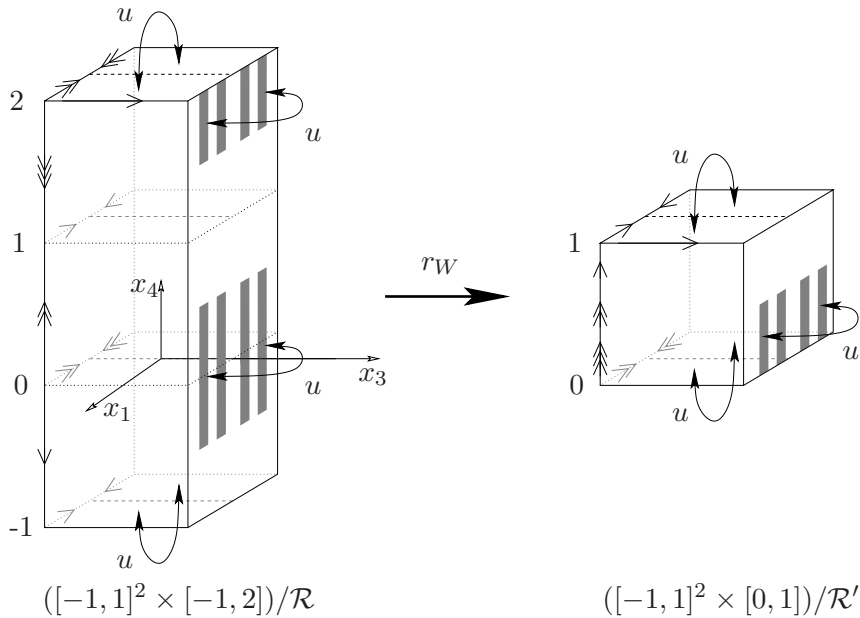} 
  \caption{$W = H^0 \cup 2 H^1$ as 3-fold branched cover
    of $\B^4$, illustrated by the 3-dimensional case (with
    coordinate directions~$x_1$,~$x_3$, and~$x_4$). The areas
    $\mathcal{A}\times([-1/2,1/2]\cup[3/2,2])$, respectively
    $\mathcal{A}\times [0,1/2]$ are shaded, and the arrows on the edges indicate the orientations of the image of an edge under~$r_W$.\label{fig:W_as_cover}}
\end{figure}

Similarly we obtain the base space of~$r_W$ from $[-1,1]^3\times[0,1]$ by
identifying a point~$x \in [-1,1]^3\times[0,1]$ with its image~$u(x)$
if~$x$ lies in the top or bottom facet, that is, if $x_4=0$ or
$x_4=1$, or if $(x_1,x_2,x_3) \in \mathcal{A}$ and
$x_4\in[0,\frac{1}{2}]$. Hence we have
$([-1,1]^3\times[0,1])/\mathcal{R}'\cong\B^4$ , where~$\mathcal{R}'$ is the
equivalence relation given by the identifications described above.
Figure~\ref{fig:W_as_cover} illustrates the 3-dimensional case.

Now we are ready to define the covering map~$r_W$. For simplicity
we identify~$W$ and~$\B^4$ with $([-1,1]^3\times[-1,2])/\mathcal{R}$,
respectively $([-1,1]^3\times[0,1])/\mathcal{R}'$, and let~$[x]$
denote the equivalence class of~$x$ in the quotient spaces~$W$
and~$\B^4$, respectively.
\[
  r_W: W \to \B^4 :
  [(x_1,x_2,x_3,x_4)] \mapsto
  \begin{cases}
    [(-x_1,x_2,x_3,-x_4)] &\text{if} \quad -1 \leq x_4 \leq 0\\
    [(-x_1,x_2,x_3,x_4)] &\text{if} \quad 0 < x_4 \leq 1\\
    [(-x_1,x_2,x_3,2-x_4)] &\text{if} \quad 1 < x_4 \leq 2
  \end{cases}
\]
The covering map~$r_W$ is well defined since it is compatible
with~$\mathcal{R}$ and~$\mathcal{R}'$, and~$r_W$ is a branched cover.
Note that the third sheet $([-1,1]^3\times[1,2])/\mathcal{R}$ is
homeomorphic to~$\B^4$ and does not contribute to the construction of
the 1-handles; it is trivial so far. Yet it will be needed in the
process of attaching the 2-handles.

We will distinguish the connected components of the branching
set of~$r_W$ as follows.  The branching set consists
of~$\lambda+1$ pairwise disjoint unknotted 2-balls~$\{P_i\}_{0\leq i\leq\lambda}$, and a  
single unknotted 2-ball~$Q$ disjoint to any of the~$P_i$.
We denote the $\lambda+1$ disjoint unknotted 2-balls by~$\mathcal{P}$, and they are given by
\[
\mathcal{P} =  P_0 \cup P_1 \cup \dots \cup P_{\lambda} = 
(\{0\}\times[-1,1]^2\times\{0\}) \cup ((\mathcal{A} /
  \mathcal{R}') \times\{0\}).
\]
The single unknotted 2-ball~$Q$ is given by
\[
Q = \{0\}\times[-1,1]^2\times\{1\}.
\]
The 2- balls $\mathcal{P} \cup Q$ 
intersect the boundary of~$\B^4$ in a system of~$\lambda+2$
unknotted and unlinked 1-spheres.

The preimage of a meridial loop $m \subset \B^4\setminus(\mathcal{P}
\cup Q)$ passing around any \mbox{$P \in \mathcal{P}$} lies in the
first and second sheet of~$r_W$, that is, $r_W^{-1}(m)$ is contained
in \mbox{$([-1,1]^3\times[-1,1])/\mathcal{R}$.} On the other hand, if
a meridial loop $m' \subset \B^4\setminus(\mathcal{P} \cup Q)$ passes
around~$Q$ we have $r_W^{-1}(m') \subset
([-1,1]^3\times[0,2])/\mathcal{R}$, and the preimage of~$m'$ lies in
the second and third sheet of~$r_W$.  Therefore the monodromy group
$\M_{r_W}$ of~$r_W$ is isomorphic to the symmetric group~$\Sigma_3$ on
three elements.  In the following we label the sheets~0,~1, and~2, and
we assume~$m$ and~$m'$ to correspond to the (generating)
transpositions $(0,1)\in \Sigma_3$ and $(0,2)\in \Sigma_3$,
respectively.

\subsubsection*{Attaching 2-handles.}
Note the (non-trivial) fact that the 2-handles $\{H^2_i\}_{1\leq i\leq
  \mu}$ may be attached independently to $W = H^0 \cup \lambda H^1$.
Hence we may assume $M_A=W\cup_h H^2$ to be obtained from~$W$ by
adding a single 2-handle $H^2\cong\B^4$.  The 2-handle~$H^2$ is
attached to~$W$ along a solid 3-torus $\Sph^1 \times \B^2$. To be more
precise, a solid 3-torus $\Sph^1 \times \B^2 \subset \partial H^2$ is
embedded into~$\partial W$ via the attaching map
$h:\Sph^1 \times \B^2\to\partial W$.  The attaching map~$h$ is
determined by the image of the meridian $\Sph^1 \times \{0\} \subset
\Sph^1 \times \B^2$, a knot~$L$ in~$\partial W$. Using isotopy the
knot~$L$ may be placed in~$\partial W$ such that its image $A = r_W(L)
\subset
\partial \B^4$ under~$r_W$ is an arc which intersects the branching
set $\mathcal{P} \cup Q$ of~$r_W$ as follows: The arc~$A$ intersects
the $\lambda + 1$ connected components $\mathcal{P}$ in its end points
only and does not intersect~$Q$ at all.  Conversely, the preimage
$r_W^{-1}(A)$ of~$A$ is the knot~$L$ and a disjoint arc~$A'$.  The
restriction~$r_W\!\left|_L\right.$ is a 2-fold branched cover of~$A$,
and~$r_W\!\left|_{A'}\right.$ is a homeomorphism corresponding to the
third, trivial sheet of~$r_W$.

In order to represent $M_A$ as a 3-fold branched
cover $r_A:M_A\to\B^4$, we attach another 4-ball~$D$
to~$\B^4$ along the 3-dimensional neighborhood $r_W\!\circ h (\Sph^1
\times \B^2)$ of~$A \subset \partial \B^4$.  This neighborhood of~$A$
is homeomorphic to~$\B^3$ if the domain $\Sph^1 \times \B^2 \subset
\partial H^2$ of~$h$ is chosen sufficiently small, and the resulting
base space remains homeomorphic to~$\B^4$. 
The preimage of~$D$ is a collection of
three copies of~$D$, two of which form the 2-handle~$H^2$ attached
to~$W$ via~$h$. The third copy~$D'\cong\B^4$ is attached to~$W$
along a 3-dimensional neighborhood of~$A'$,
that is, we attach~$D'$ to~$W$ along a 3-ball, and attaching~$D'$ does
not alter the homeomorphic type of~$M_A$.

\begin{figure}[t]
  \center
  \psfrag{P1}{$P_0$}
  \psfrag{P2}{$P_1$}
  \psfrag{P3}{$P_2$}
  \psfrag{R}{$Q$}
  \psfrag{A1}{\hspace{-1cm}$A_1\subset\overline{A_1}$}
  \psfrag{A2}{$\overline{A_2}$}
  \includegraphics[width=8.5cm]{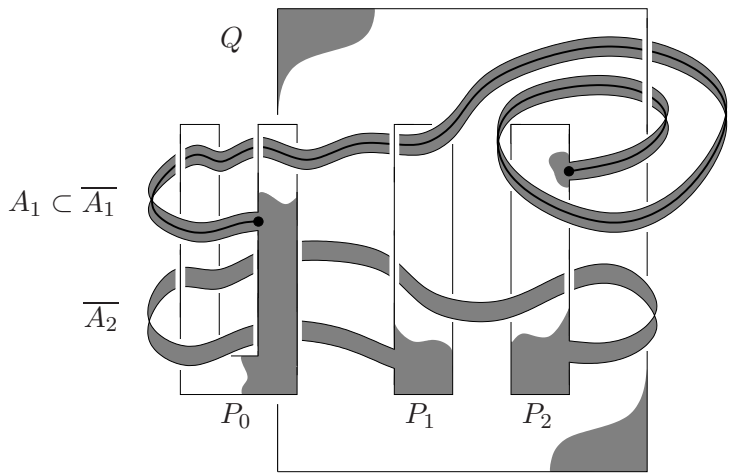} 
  \caption{Immersion of a ribbon manifold with two ribbons
    $\overline{A_1}$ and~$\overline{A_2}$. Additionally the arc $A_1\subset\overline{A_1}$
    is pictured.\label{ribbon_manifold}}
\end{figure}

The resulting branching set of $r_A$
is the union of the branching set
of~$r_W$, the 2-balls $\mathcal{P} \cup Q$, and a 2-ball
$\overline{A} \supset A$ attached to~$\mathcal{P}$ along two arcs~$a$
and~$a'$ in the boundary of~$\mathcal{P}$, a
\emph{ribbon manifold}; see Figure~\ref{ribbon_manifold}.
The two arcs~$a$ and~$a'$ are neighborhoods in~$\partial\mathcal{P}$
of the two endpoints $A\cap\mathcal{P}$ of~$A$. Note that~$r_A$
is a ``proper'' 3-fold branched cover (the
third sheet is non-trivial), since although~$\overline{A}$ does not
intersect~$Q$, it might ``weave around''~$Q$ (and in fact also around
\mbox{any~$P \in \mathcal{P}$}). 

Fix a set of meridial loops as generators of
$\pi_1(\B^4\setminus(\mathcal{P}\cup Q),y_0)$, that is, choose one
meridial loop around each of the 2-balls in~$\mathcal{P}$, and one
meridial loop around~$Q$.  Let~$P,P'\in\mathcal{P}$ with $a\subset P$
and $a'\subset P'$, and let
$\beta,\beta'\in\pi_1(\B^4\setminus(\mathcal{P}\cup Q),y_0)$ be the
generators corresponding to the meridial loops around~$P$ and~$P'$,
respectively. Then adding the ribbon~$\overline{A}$ to the branching
set introduces a new relation to the fundamental group, that is,
the group $\pi_1(\B^4\setminus(\mathcal{P}\cup Q\cup\overline{A}),y_0)$ differs
from $\pi_1(\B^4\setminus(\mathcal{P}\cup Q),y_0)$ by the relation
\begin{equation}\label{equ:ribbon_relation}
  \beta\alpha=\beta',
\end{equation}
where the element $\alpha\in\pi_1(\B^4\setminus(\mathcal{P}\cup Q),y_0)$ corresponds to the
way~$\overline{A}$ weaves around $\mathcal{P}\cup Q$.
We summarize the construction above by the following Theorem~\ref{thm:montesinos_1}.

\begin{thm}[{Montesinos~\cite[Theorem~6]{montesinos:CSR}}]
  \label{thm:montesinos_1}
  Each $4$-manifold $M_A = H^0 \cup \lambda H^1 \cup \mu H^2$ is a
  simple $3$-fold branched cover of the $4$-ball, the branching set being
  a ribbon manifold.
\end{thm}

\subsubsection*{Construction of $H$.}
The construction of the cobordism~$r_H: H \to \Sph^3\times [0,1]$
is rather straight forward once we have established its existence, which
is provided by the Theorems~\ref{thm:montesinos_2}
and~\ref{thm:piergallini_1}. Note that
the branched cover~$r_H: H \to \Sph^3\times [0,1]$ is already
defined on the boundary of~$H$ by the restrictions given in
Equation~\eqref{equ:cobordism}: The boundary of~$H$ is the disjoint
union of two copies of the 3-manifold $\g \; \sharp \; (\Sph^1 \times
\Sph^2)$, and the branching sets of the restrictions
$r_A\left|_{\partial M_A}\right.$ and $r_B\left|_{\partial
    M_B}\right.$ are two links~$L_A$ and~$L_B$, respectively.

\begin{figure}[t]
  \center
  \psfrag{C+}{$C^+$}
  \psfrag{C-}{$C^-$}
  \psfrag{N+}{$N^{+}$}
  \psfrag{N-}{$N^{-}$}
  \psfrag{(i,j)}{\small$(i,j)$}
  \psfrag{(k,l)}{\small$(k,l)$}
  \psfrag{(i,k)}{\small$(i,k)$}
  \psfrag{(j,k)}{\small$(j,k)$}
  \includegraphics[width=10.3cm]{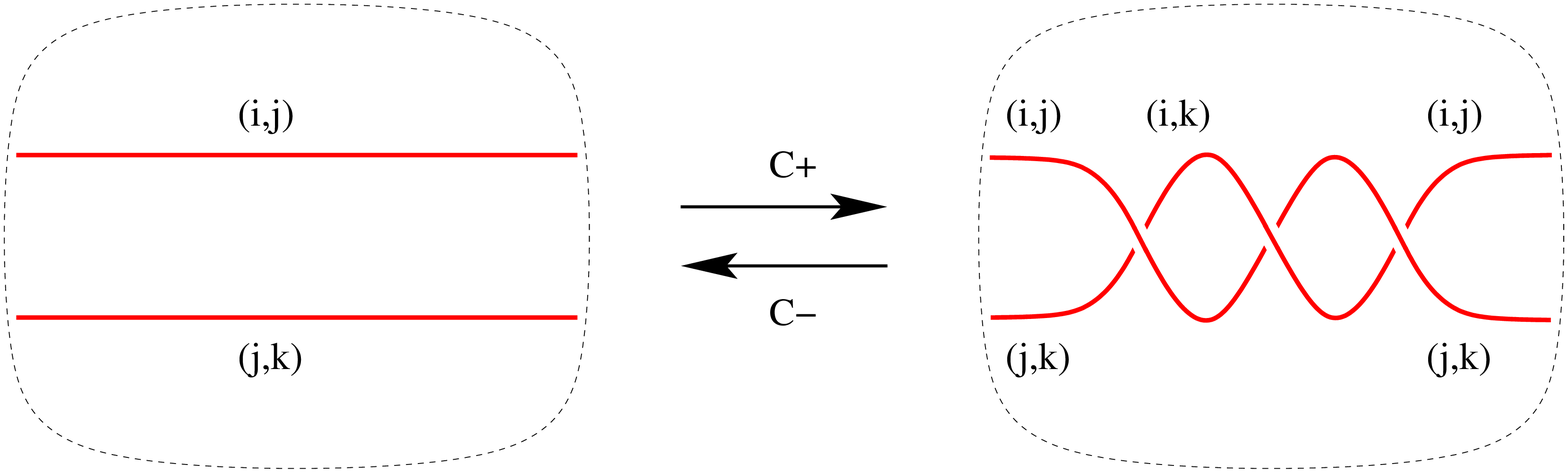}
  \includegraphics[width=10.3cm]{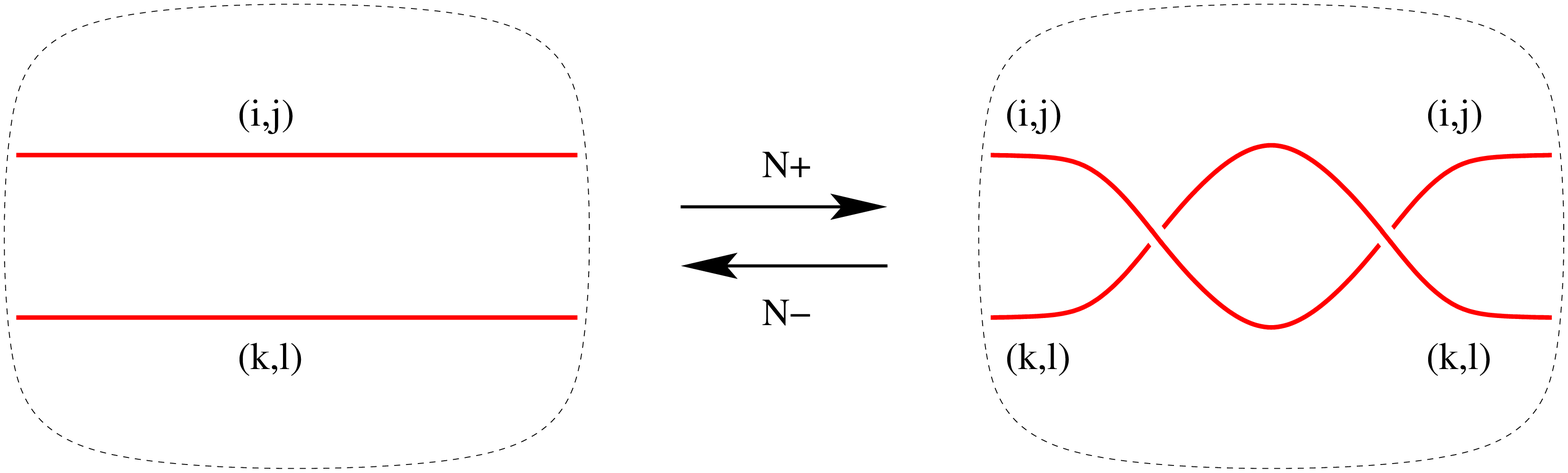}
  \caption{The moves~$C^\pm$ and~$N^\pm$.\label{fig:CN_moves}}
\end{figure}

In general, any closed oriented 3-manifold~$W^3$ arises as a simple 3-fold branched
cover of~$\Sph^3$ branched over a link~$L$, and the monodromy
group~$\M$ of the branched cover is isomorphic to a subgroup
of~$\Sigma_3$ (generated by transpositions);
see~\cite{hilden:TFBC,montesinos:MFBC}.  After adding a fourth trivial
sheet and thus a new generating transposition to~$\M$,~$\M$ becomes isomorphic
to a subgroup of~$\Sigma_4$.  Consider a generic projection of~$L$ to
the plane with marked over and under crossings. Such a projection is
called a \emph{diagram} of~$L$, and we call a strand which is not
crossed by other strands of the diagram a \emph{bridge}.  Fix a set of
meridial loops around the bridges of the diagram as generators of
$\pi_1(\Sph^3\setminus L)$, and we identify the meridial loops around the
bridges with transpositions in~$\M$ via the monodromy homomorphism
$\m:\pi_1(\Sph^3\setminus L)\to\M$. Hence we can think of~$L$ as a
\emph{colored diagram}: A bridge~$b$ of the diagram is colored~$(i,j)$ if
the meridial loop around~$b$ corresponds to the
transposition~$(i,j)\in\Sigma_4$.  Further we define the moves~$C^\pm$
and~$N^\pm$ on a colored link as in Figure~\ref{fig:CN_moves}.

\begin{thm}[{Montesinos~\cite[p.~345]{montesinos:MIC}}]
  \label{thm:montesinos_2}
  Let $p_1, p_2: W^3 \to \Sph^3$ be $4$-fold branched covers (coming
  from $3$-fold covers by the addition of a trivial sheet) such that it
  is possible to pass from the branching set~$L_1$ of~$p_1$ to the
  branching set~$L_2$ of~$p_2$ by a sequence of moves~$C^\pm$
  and~$N^\pm$. Then~$p_1$ and~$p_2$ are cobordant and the branching
  set of the cobordism is an embedded PL $2$-manifold with a cusp
  singularity (a cone over the trefoil) for each~$C^\pm$-move and a
  node singularity (a cone over the Hopf link) for each~$N^\pm$-move.
\end{thm}


To understand the main idea of the proof it suffices to look at two
branched covers $p_1, p_2: W^3 \to \Sph^3$ such that their branching
sets~$L_1$ and~$L_2$ differ by exactly one~$C^\pm$- or~$N^\pm$-move
$m$. Let~$U\subset\Sph^3$ be a closed neighborhood of the
move~$m$, that is, $L_1\setminus U\subset \Sph^3\setminus U$ and
$L_2\setminus U\subset \Sph^3\setminus U$ are equivalent, and replacing
$L_1\cap U$ by $L_2 \cap U$ realizes the move~$m$.
The branching set in $(\Sph^3\setminus U) \times [0,1]$ is
$(L_1\setminus U) \times [0,1]\cong(L_2\setminus U) \times [0,1]$. If~$m$ is a~$C^\pm$-move then the
intersection of the branching set $(L_1\setminus U) \times [0,1]$ with
the boundary of $U\times [0,1]$ is the trefoil, otherwise the
intersection is the Hopf link. In order to complete the base space of
our cobordism,
we replace $U\times [0,1]$ by a 4-ball with a cone over
the trefoil or the Hopf link, respectively, as a branching set.

Theorem~\ref{thm:montesinos_2} together with the following
Theorem~\ref{thm:piergallini_1} establish the existence of the
cobordism~$r_H$, and completes the construction of the branched
cover $r:M\to\Sph^4$. As observed by
Montesinos~\cite{montesinos:MIC}, Theorem~\ref{thm:piergallini_2} then
follows immediately.

\begin{thm}[{Piergallini~\cite[Theorem~A]{piergallini:FMS4}}]
  \label{thm:piergallini_1}
  Any two branching sets of $4$-fold branched covers $p_1, p_2: W^3 \to \Sph^3$ obtained
  from $3$-fold branched covers by adding a fourth, trivial sheet, which
  represent the same $3$-manifold~$W^3$, are related by a finite
  sequence of moves~$C^\pm$ and~$N^\pm$.
\end{thm}

The proof extends over two papers by Piergallini. In~\cite{piergallini:CM} the
number of different moves needed to relate any two such branching sets
via a finite sequence of moves
is brought down to four. Then in~\cite{piergallini:FMS4} each of these four
moves is realized by a finite sequence of $C^\pm$-
and~$N^\pm$-moves, and the usage of a fourth, trivial sheet, thus establishing
Theorem~\ref{thm:piergallini_1}.

\begin{thm} \label{thm:piergallini_2}
  Every closed oriented PL
  $4$-manifold is a simple $4$-fold branched cover of the $4$-sphere branched
  over a immersed PL-surface with a finite number of cusp and node
  singularities.
\end{thm}

\subsection{4-Manifolds as partial unfoldings}
\label{sec:4mf_as_unfoldings}

Let~$M$ be a closed oriented PL 4-manifold, and let $r:M\to\Sph^4$
be the 4-fold branched cover with branching set~$F$ described in
Section~\ref{sec:4mf_as_branched_covers}. Hence~$F$ is an immersed
PL-surface with a finite number of cusp and node singularities by
Theorem~\ref{thm:piergallini_2}. In Theorem~\ref{thm:ch_4mf_main} we
construct a triangulation~$S$ of~$\Sph^4$ such that the branched cover
given by the projection $p:\widehat{S}\to S$ is equivalent to $r:M\to\Sph^4$. In
particular,~$\widehat{S}$ is PL-homeomorphic to~$M$.  Recall that we
refer to the (unique non-trivial) connected component of the partial
unfolding PL-homeomorphic to~$M$ as \emph{the} partial unfolding.

We outline the construction of~$S$. The branched cover~$r$ is
characterized by~$F$ and the monodromy isomorphism $\m_r
:\pi_1(\Sph^4\setminus F,y_0)\to\Sym(r^{-1}(y_0))$, where~$y_0$ is a
point in $\Sph^4\setminus F$; see Section~\ref{sec:bc} and
Theorem~\ref{thm:bc}. Therefore we construct~$S$ such that there is a
homeomorphism of pairs~$\varphi:(\Sph^4,F)\to (|S|,|S_{\odd}|)$ and~$\varphi$ induces a group
isomorphism $\varphi_*:\pi_1(\Sph^4\setminus
F,y_0)\to\pi_1(|S|\setminus|S_{\odd}|,\varphi(y_0))$. Further, assume
that $\varphi(y_0)$ is the barycenter of some facet $\s_0\in S$. We
construct~$S$ such that the following
Diagram~\eqref{equ:monodr_and_proj} commutes for some bijection
$\iota:r^{-1}(y_0)\to V(\s_0)$ and the induced group isomorphism
$\iota_*:\M_r\to\Pi(S,\s_0)$.
\begin{equation}\label{equ:monodr_and_proj}
\xymatrix{
  \pi_1(\Sph^4\setminus F,y_0) \ar[r]^{\varphi_*} \ar[d]^{\m_r} 
  & \pi_1(|S|\setminus|S_{\odd}|,\varphi(y_0)) \ar[d]^{\h_S}\\
  \M_r \ar[r]^{\iota_*} & \Pi(S,\s_0)\\
}
\end{equation}

This establishes Theorem~\ref{thm:ch_4mf_main}, since the partial unfolding of a
nice simplicial complex
is a branched cover by Corollary~\ref{cor:unf_are_covers}, and since
\mbox{$\varphi(F)=|S_{\odd}|$} and commutativity of
Diagram~\eqref{equ:monodr_and_proj} ensures that~$p:\widehat{S}\to S$ and~$r:M\to\Sph^4$
are indeed equivalent by Theorem~\ref{thm:partial_unf}. The PL-properties follow once we proved~$S$ to be a combinatorial manifold.

The construction of~$S$ follows closely the construction of the
branched cover $r:M\cong M_A\cup H\cup M_B \to \B^4\cup(\Sph^3\times
[0,1])\cup\B^4$ reviewed in Section~\ref{sec:4mf_as_branched_covers}:
First the combinatorial 4-balls~$D_A$ and~$D_B$ are constructed such
that $\widehat{D}_A\cong M_A$ and $\widehat{D}_B\cong M_B$,
respectively.  The resulting complex~$T_1$ is the disjoint union
of~$D_A$ and~$D_B$. For each $C^\pm$- and $N^\pm$-move~$m$ needed to
relate the odd subcomplexes of $\partial D_A$ and $\partial D_B$ we
then attach a 4-ball~$D_m$ to~$D_A$ such that the partial unfolding of
$D_A\cup D_H=D_A\cup(\bigcup_m D_m)$ is PL-homeomorphic to $M_A \cup H$. We refer to
the simplicial complex constructed as~$T_2= D_A\cup D_H \dotcup D_B$, and we have $T_1\subset
T_2$. In a last step we triangulate the
remaining space $\Sph^4\setminus |T_2|\cong\Sph^3\times [0,1]$,
attaching~$D_B$ to $D_A\cup D_H$.  This yields $T_3=S$. In
each step~$T_1$,~$T_2$, and~$T_3$ of the construction of~$S$ we have
to ensure
\begin{enumerate}[(I)]
\item that $\varphi(F)\cup|T_i|=|\odd(T_i)|$ and
\item that Diagram~\eqref{equ:monodr_and_proj} restricted to~$T_i$ commutes.
\end{enumerate}

Note that each of the complexes~$T_i$ has to be nice for~$\h_{T_i}$ to be well defined.
Finally we may assume~$T_i$ to be a
sufficiently fine triangulation. A fine triangulation can be obtained
by anti-prismatic subdivision
of faces at any stage of the
construction by Proposition~\ref{prop:anti-prismatic}.

\subsubsection*{Construction of $T_1=D_A\dotcup D_B$.}
We begin with constructing a triangulation~$D_W$ of the base space of
the branched cover $r_W: W = H^0 \cup \lambda H^1 \to \B^4$, that
is, $\widehat{D_W}\cong W$. Then we modify~$\odd(D_W)$ by adding
the branching set which produces the~$\mu$ 2-handles in
order to construct a triangulation~$D_A$ of the base space of $r_A:
M_A = H^0 \cup \lambda H^1 \cup \mu H^2\to \B_4$, that is $\widehat{D_A}\cong M_A$. To this end let~$C$
be a sufficiently fine triangulated foldable combinatorial 4-ball obtained
via the iterated barycentric subdivision
of a 4-simplex. Since~$C$ arises
as a barycentric subdivision there is a natural 5-coloring of the
vertices of~$C$ by coloring each vertex~$v\in C$ by the dimension of
the original face subdivided by~$v$.
Therefore~$\partial C$ lies in the
$\{0,1,2,3\}$-skeleton, and
vertices colored~4 appear only in the interior of~$C$. The
triangulation~$D_W$ of~$\B^4$ (and later the triangulation~$D_A$) is
obtained from~$C$ by a series of stellar subdivisions of edges.
To cut down on notation we keep referring to our complex by~$C$
throughout all stages of the construction, and~$C$ is 6-colorable
assigning a new color to all new vertices while preserving the
original coloring otherwise; see Proposition~\ref{prop:odd_subcomplex}.

In order to specify the isomorphism
$\iota_*:\M_r\to\Pi(S,\s_0)$ in
Equation~\eqref{equ:monodr_and_proj} fix a facet~$\s_0\in C$ and
let~$\iota$ map the element $x_i\in p^{-1}(y_0)$ contained in the
$i$-th sheet of~$r$ to the vertex of~$\s_0$ colored
$j\in\{0,\dots,4\}$ via the permutation
\[
\left(
\begin{matrix}
  \;0 \quad 1 \quad 2 \quad 3 \quad 4\;\\
  \;3 \quad 1 \quad 2 \quad 4 \quad 0\;
\end{matrix}
\right)
\]
We will keep~$\s_0$ fixed throughout the construction
of~$S$. Although the choice for~$\iota$ may seem arbitrary, it turns
out to be useful when applying Lemma~\ref{lem:odd_subcomplex} in the
construction of~$D_W$.

Recall that subdividing an edge~$e$ in the $\{i,j\}$-skeleton yields
$\Link(e)$ as the odd subcomplex in the complementary skeleton, that
is, in the $(\{0,\dots,4\}\setminus\{i,j\})$-skeleton; see
Proposition~\ref{prop:odd_subcomplex}.  A projectivity around a
triangle in~$\Link(e)$ exchanges the two vertices of~$\s_0$
colored~$i$ and~$j$. Via~$\iota^{-1}$ such a projectivity corresponds
to exchanging the elements of~$r^{-1}(y_0)$ contained in the
sheets of~$r$ labeled $\iota^{-1}(i)$ and $\iota^{-1}(j)$.

We first realize the 2-balls in~$\mathcal{P}$ as the odd
subcomplex in the $\{0,2,4\}$-skeleton, since they correspond
via~$\iota^{-1}$ to the
transposition~$(0,1)$ in~$\M_{r_W}$. To this end we embed for each
$P\in\mathcal{P}$ a 3-ball~$F_P$ in the $\{0,2,3,4\}$-skeleton with
$\partial F_P$ in the $\{0,2,4\}$-skeleton, and $P\cong\closure(\partial F_P \setminus
\partial C)$. Such an embedding of~$F_P$ exists by
Lemma~\ref{lem:odd_subcomplex} since we assume~$C$ to be sufficiently
finely triangulated, and we choose the $\{F_P\}_{P\in\mathcal P}$
pairwise disjoint. Now we obtain~$\mathcal{P}$ as the odd subcomplex
by stellar subdivision of $\{1,3\}$-edges following
Proposition~\ref{prop:odd_subcomplex}.

The odd subcomplex representing~$Q$ is built in a similar fashion in
the $\{0,1,4\}$-skeleton, since~$Q$ corresponds via~$\iota^{-1}$ to the
transposition~$(0,2)$ in~$\M_{r_W}$. The 3-ball~$F_Q$ with
$Q\cong$ \mbox{$\closure(\partial F_Q \setminus\partial C)$} 
is embedded in the $\{0,1,3,4\}$-skeleton with~$\partial F_Q$ in the
$\{0,1,4\}$-skeleton. Proposition~\ref{prop:odd_subcomplex} is
applicable since~$\mathcal{P}$ and~$F_Q$ are
disjoint. Now~$Q$ is realized as the odd subcomplex in the
$\{0,1,4\}$-skeleton by subdividing $\{2,3\}$-edges. 
This completes the construction of~$D_W$. The odd subcomplex
intersects~$\partial C$ in a system of $\lambda + 1$ unknotted and
unlinked~$\Sph^1$ in the $\{0,2\}$-skeleton representing~$\partial
\mathcal{P}$, and a single unknotted and unlinked~$\Sph^1$ in the
$\{0,1\}$-skeleton representing~$\partial Q$.

Finally we have to add the~$\mu$ ribbons to the odd
subcomplex in order to construct~$D_A$. To this end let~$y_0$ be the barycenter of~$\s_0$, and
fix a set of meridial loops as generators of
$\pi_1(C\setminus(\mathcal{P}\cup Q),y_0)$, that is, choose one meridial loop
around each of the 2-balls in~$\mathcal{P}$, and one meridial loop around~$Q$.
Further assume that the generators do not intersect
the collection of 3-balls $\{F_P\}_{P\in\mathcal{P}}\cup F_Q$. Then a projectivity along the image
under~$\h_C$ of a generator around a
2-ball $P\in\mathcal{P}$ exchanges the vertices colored~1 and~3
of~$\s_0$, and a projectivity along the image under~$\h_C$ of the generator around the
2-ball~$Q$ exchanges the vertices colored~2 and~3.

Now let~$A\in\partial C$ be the arc corresponding to a
ribbon~$\overline{A}$ and let~$a\subset P$ and~$a'\subset P'$ be the
intersection of~$\overline{A}$ with~$\mathcal{P}$ as described in
Section~\ref{sec:4mf_as_branched_covers}. Further let~$\beta$
and~$\beta'$ be the elements of $\pi_1(C\setminus(\mathcal{P}\cup
Q),y_0)$ corresponding to the meridial loops around~$P$ and~$P'$. In
order to apply Proposition~\ref{prop:odd_subcomplex} choose a regular
4-dimensional neighborhood~$U_A$ of~$A$ in~$C$. (Provided that~$C$ is
sufficiently fine triangulated one may choose $U_A=\bigcup_{v\in
  A}\Star_C(v)$.) The neighborhood~$U_A$ is 5-colorable since
$\odd({U_A})=\emptyset$, and we may choose the coloring such that it
coincides with the coloring of~$C$ in neighborhoods of~$a$ and~$a'$,
respectively. The later assumption holds since $\beta \alpha = \beta'$
holds by Equation~\eqref{equ:ribbon_relation},
where~$\alpha\in\pi_1(C\setminus(\mathcal{P}\cup Q),y_0)$ corresponds
to the way~$A$ weaves around $\mathcal{P}\cup Q$. Observe that the
5-coloring~$U_A$ does not coincide with the coloring of~$C$ in
general. It changes corresponding to the way~$A$ weaves around the
2-balls $\mathcal{P}\cup Q$.

Now choose a 3-ball~$F_{\overline{A}}$ according to
Proposition~\ref{prop:odd_subcomplex} in the $\{0,2,3,4\}$-skeleton of~$U_A$
with $\partial F_{\overline{A}}$ in the $\{0,2,4\}$-skeleton, and
$\overline{A}\cong\closure(\partial F_{\overline{A}}\setminus\partial C)$.
If we color the vertices of~$U_A$
by the coloring of~$C$, then in general~$\partial F_{\overline{A}}$ is partly embedded in the
$\{0,2,4\}$-, $\{0,1,4\}$-, and $\{0,3,4\}$-skeleton, reflecting the fact that different parts
of the ribbon correspond to different transposition~$(0,1)$,~$(0,2)$, and~$(1,2)$. The intersection
of~$\overline{A}$ with~$\mathcal{P}$ however is always contained in
the $\{0,2\}$-skeleton.

The ribbon~$\overline{A}$ is added to the odd subcomplex by stellar subdividing edges in the $\{1,3\}$-skeleton of~$U_A$ by Proposition~\ref{prop:odd_subcomplex}.
Adding all ribbons $\bigcup_{i=1}^\mu\overline{A}_i$ to the odd
subcomplex completes the construction of~$D_A$.

The simplicial 4-balls~$D_A$ and~$D_B$ are indeed combinatorial
4-balls (and hence nice) since they are constructed by subdivision of faces from the
4-simplex, and they meet conditions~(I)
and~(II)
by construction. We have $T_1=D_A\dotcup
D_B$, and we summarize the construction of~$D_A$ by the following proposition.

\begin{prop}\label{prop:M_A_as_unfolding}
  For each PL $4$-manifold $M_A = H^0 \cup \lambda H^1 \cup \mu H^2$ there is a
  combinatorial $4$-ball~$D_A$ such that one of the connected components of the
  partial unfolding~$\widehat{D_A}$ is PL-homeomorphic to~$M_A$. The
  projection $\widehat{D_A}\to D_A$ is a simple
  $3$-fold branched cover with a ribbon manifold as a branching set.
\end{prop}

\subsubsection*{Construction of $T_2=D_A\cup D_H\dotcup D_B$.}
For the construction of the cobordism $r_H \to
\Sph^3\times [0,1]$ we need~$r_A$ and~$r_B$ to be 4-fold branched
covers obtained from a 3-fold branched cover by adding a trivial
sheet.  The fourth sheet is obtained by adding a 2-ball in the
$\{0,1,2\}$-skeleton to the odd subcomplex via stellarly subdividing
edges in the $\{3,4\}$-skeleton by
Proposition~\ref{prop:odd_subcomplex} and
Lemma~\ref{lem:odd_subcomplex}.  A projectivity along a closed
facet path based at~$\s_0$ around a triangle of the
newly added odd subcomplex exchanges the vertices of~$\s_0$ colored~3
and~4, and corresponds via~$\iota^{-1}$ to the transposition~$(0,3)$
in~$\M_r$.

We first construct $D_A\cup D_H$ such that its
partial unfolding yields $M_A\cup H$. In particular the odd subcomplex
of the boundary of $D_A\cup D_H$ is equivalent to
the odd subcomplex of~$D_B$. To this end a combinatorial 4-ball~$D_m$ is attached
to~$\partial D_A$ successively for
each of the~$C^\pm$- and $N^\pm$-moves required to relate the odd
subcomplex of~$\partial D_A$ and~$\partial D_B$. The 4-ball~$D_m$
realizes~$m$ in the sense that the odd subcomplexes of
$\partial D_A$ and $\partial (D_A\cup D_m)$ differ by the move~$m$.
This produces the
triangulation~$T_2$. We then ``identify'' the
boundaries of~$D_B$ and $D_A\cup(\partial D_A\times [0,1])$, thus
completing the triangulation $S=T_3$. Keep in mind that we have to
ensure conditions~(I) and~(II)
to be valid throughout the construction.

The combinatorial 4-ball~$D_m=\Cone(S_m)$ is constructed as the cone over a combinatorial
3-sphere~$S_m$ with a trefoil knot or Hopf link as (colored!)
odd subcomplex, respectively.
The resulting odd
subcomplex $\odd(D_m)$ is a cusp or a node singularity depending on
whether~$m$ is a~$C^\pm$- or~$N^\pm$-move, since $\odd(D_m)=\odd(\Cone(S_m))=\Cone(\odd(S_m))$ holds.

\begin{figure}[t]
  \centering
  \begin{overpic}[width=.43\textwidth]{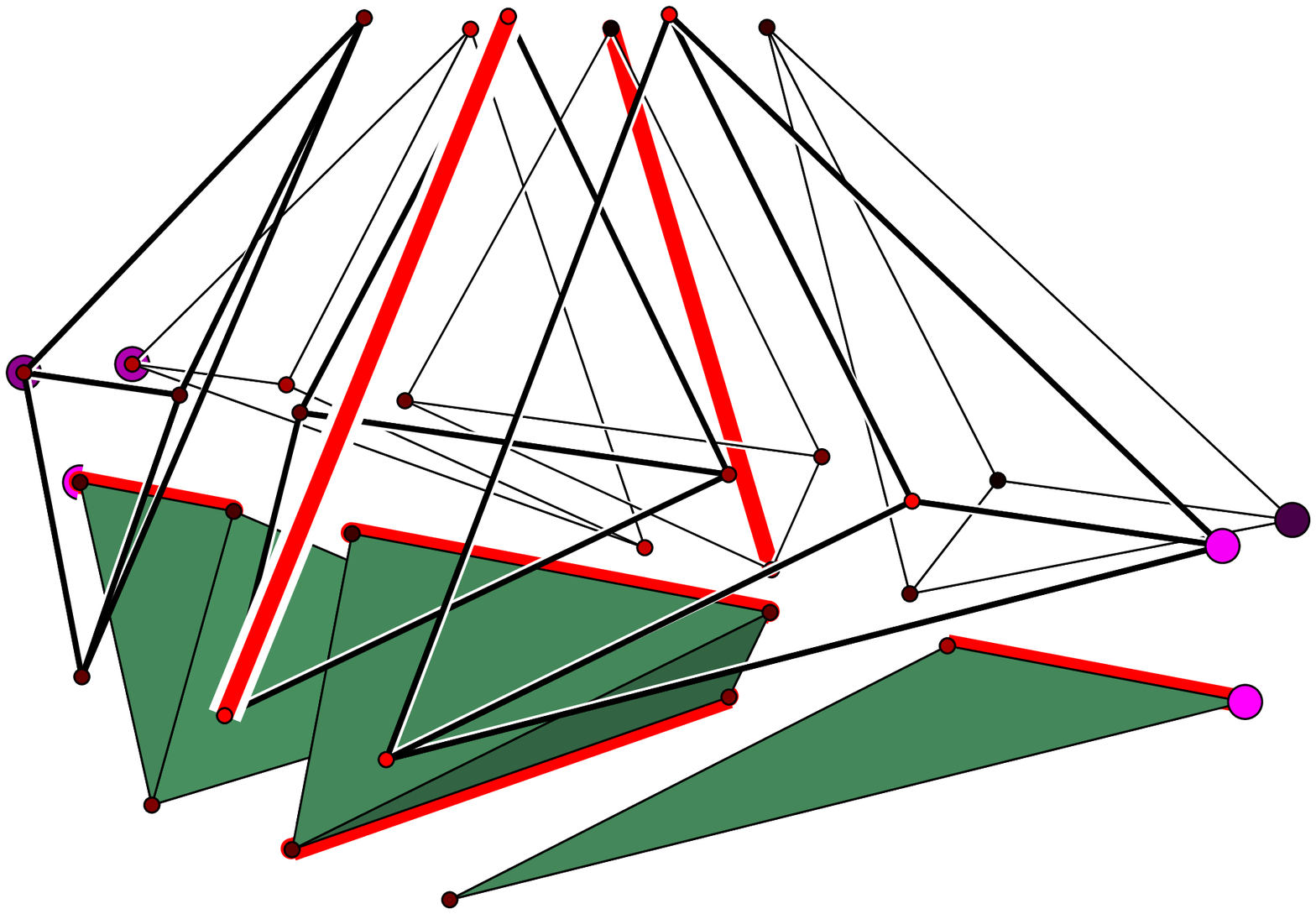}
    \put(-3,40){\textbf{\textit{A}}}
    \put(99,23){\textbf{\textit{A}}}
  \end{overpic}
  \includegraphics[width=.42\textwidth]{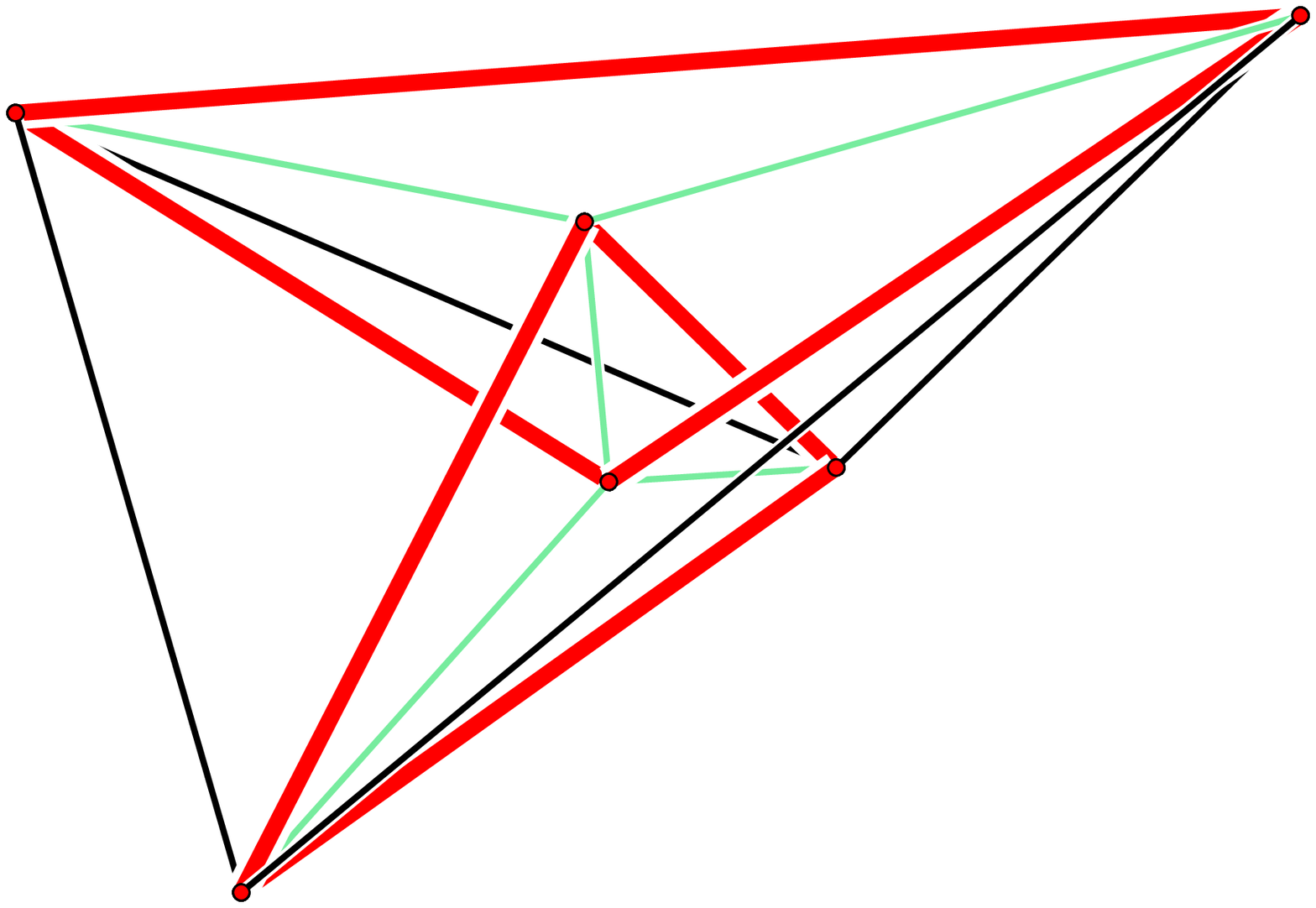}
  \caption{Construction of~$S_m$ with the Hopf link as odd subcomplex and
    $\Pi(S_m)\cong\Sigma_2\times\Sigma_2\subgr \Sigma_4$ as the cone
    (with apex~$A$) over the triangulated bipyramid, and as the boundary complex of the direct-sum of two triangles, pictured as its Schlegel diagram.
    The odd subcomplex is marked.
    \label{fig:Sm_hopf}}
\end{figure}

In general, the sphere~$S_m$ may be obtained
following the construction by Izmestiev~\& Jos\-wig~\cite{izmestiev_joswig:BC}.
Alternatively, an explicit triangulation of~$S_m$ for a $C^\pm$-move~$m$ is available as 
an electronic model (\texttt{polymake}~\cite{polymake} file) by~\cite{witte:DIPL}.
For a $N^\pm$-move~$m$ a triangulation of~$S_m$ may be obtained as the cone
over the triangulated bipyramid, respectively as the boundary complex of 
the direct sum of two triangles; see Figure~\ref{fig:Sm_hopf}.

\begin{rem}\label{rem:unf_commutes}
  Observe that the operations ``coning'' and ``partially unfolding'' commute and
  that for~$m$ being either a $C^\pm$- or $N^\pm$-move the partial unfolding~$\widehat{S_m}$ is again
  a combinatorial 3-sphere.
\end{rem}

It remains to show how to attach~$D_m$ to~$D_A$. Choose 3-dimensional neighborhoods
$U\subset\partial D_A$ and $U'\subset\partial D_m=S_m$, such that
replacing~$U$ by $S_m\setminus U'$ realizes the
move~$m$.
Now the move~$m$ is realized by identifying~$|U|$ and~$|U'|$.
Since the triangulations~$U$ and~$U'$ are non-equal in general, we
triangulate the space $|U|\times[0,1]$, such that~$U$ triangulates
$|U|\times\{0\}$ and~$U'$ triangulates $|U|\times\{1\}$, and such that
the odd subcomplex is equivalent to the prism over~$|U_{\odd}|$.

Attaching~$D_m$ to~$D_A$ by
identifying~$|U|$ to~$|U'|$ is similar to the last remaining step in the
construction of~$S$, where~$M_B$ is attached to~$M_A\cup H$
via identifying~$|\partial (M_A\cup
H)|$ and~$|\partial M_B|$. We explain
how to realize the identification of~$|U|$ and~$|U'|$, respectively
of~$|\partial (M_A\cup H)|$ and~$|\partial M_B|$, via extending the
triangulations~$U$
and~$U'$, respectively~$\partial M_B$
and~$\partial (M_A\cup H)$ in a more general setting,
thus completing the construction of~$S$.

\subsubsection*{Attaching along color equivalent subcomplexes.}
Consider two combinatorial 4-manifolds~$K$ and~$K'$ and combinatorial 3-manifolds (possibly
with boundary)~$U\subset\partial K$ and~$U'\subset\partial K'$
with~$U\cong U'$.  Assume that there are color equivalent regular
4-dimensional neighborhoods~$N$ and~$N'$ of~$U$, respectively~$U'$, such
that~$|N_{\odd}|$ (and hence also~$|N'_{\odd}|$) is equivalent to
$|U\cup N_{\odd}|\times[0,1]$, and~$N_{\odd}$ is a locally flat
combinatorial $2$-manifold. Note that $U_{\odd}= U\cap N_{\odd}$
does not hold in general. Further let~$\varphi: |N|\to |N'|$, $\s_0\in
N$, $\s'_0\in N'$, and $\psi:V(\s_0)\to V(\s'_0)$ as in
Equation~\eqref{equ:color_equivalent}, defining the color equivalence
of~$N$ and~$N'$.

\begin{prop}\label{prop:attaching}
  There is a triangulation~$T$ of $|U|\times[0,1]$ with $|T_{\odd}|$
  equivalent to
  $|U\cap N_{\odd}|\times[0,1]$, such that~$T$ equals~$U$ on $|U|\times\{0\}$
  and~$U'$ on $|U|\times\{1\}$, and such that
  $\odd(K\cup T\cup K')=K_{\odd}\cup T_{\odd} \cup
  K'_{\odd}$, thus in effect attaching~$K'$ to~$K$ via identification
  of~$U$ and~$U'$.  Here $K\cup T\cup K'$ denotes the union
  of~$K$,~$K'$, and~$T$, attaching~$T$ to~$K$ and~$K'$ along~$U$,
  respectively~$U'$.
  The simplicial complex $K\cup T\cup K'$ is a combinatorial $4$-manifold.
\end{prop}

In order to make the proof digestible it is split into the three
Lemmas~\ref{lem:attaching1},~\ref{lem:attaching2}, and~\ref{lem:attaching3}.
We denote a face $f\in N$ which intersects~$U$ in
all except one vertex by $f=\{g,x_g\}$, where~$g$ is a
face of~$U$ and~$x_g$ the one remaining vertex.  Faces of~$N'$
intersecting~$U'$ in all except one vertex are denoted
similarly. Throughout,~$\tau\in U$ will be a facet of~$U$, that is, a
ridge of~$N$.
  
After possible refinements of~$N$ and~$N'$ via anti-prismatic
subdivision there is a simplicial approximation~$\varphi':N\to N'$
of~$\varphi$ which does not degenerate~$\s_0$. Note that any
simplicial approximation of~$\varphi$
maps~$N_{\odd}$ to~$N'_{\odd}$, and~$U$ to~$U'$;
see~\cite[Lemma 14.4, Theorem 16.1]{munkres:EAT}.
Let~$\s\in N$ be a facet,~$\g$ a facets path in~$n$ from~$\s_0$
to~$\s$, and let~$\g'$ be the facet path in~$N'$ defined by the
non-degenerated images of facets in~$\g$. Let~$\kappa_\s$
be the last facet of~$\g'$, hence $\s'=\varphi'(\s)\subset\kappa_\s$ in general,
and $\s'=\kappa_\s$ if~$\varphi'$ does not degenerate~$\s$. We
define the bijective map $\psi_\s:V(\s)\to V(\kappa_\s)$ by
\[
\psi_\s=\langle\g'\rangle\circ\psi\circ\langle\g\rangle^{-1}.
\]
Since~$N$ and~$N'$ are color equivalent,~$\psi_\s$ is independent of
the choice of~$\g$ and hence well defined.  Further note
that $\psi_\s^{-1}\left|_{\s'}\right.$ is injective, and
that \mbox{$\psi_\s(\s\cap N_{\odd})\subset N'_{\odd}$} since~$N$ and~$N'$
are color equivalent.

Consider the following regular cell decomposition of $|U|\times
[0,1]$. First the $i$-faces of~$U$ and~$U'$ form closed $i$-cells in the
natural way.  In particular, the vertices of~$U$
and~$U'$ are the 0-cells. Now we add a closed $(i+1)$-cell~$C_f^{i+1}$
for each $i$-face~$f\in U$.
The $(i+1)$-cell~$C_f^{i+1}$ is attached to the union of all $i$-cells
along the cell decomposition of~$\Sph^i$ given by the cells~$f$ (and its proper faces),~$\varphi'(f)$
(and its proper faces), and all cells~$C_g^{j+1}$ with~$g\subset
f$ is a $j$-face.
The top dimensional cells are the 4-cells
$\{C_\tau^4\}_{\tau\in U}$ corresponding to facets of~$U$. Any two cells~$C_f^{i+1}$ and~$C_g^{j+1}$
intersect properly, that is, in the common cell corresponding to $f\cap g$, and the union of all
cells \mbox{equals~$|U|\times [0,1]$.}

We describe how to triangulate~$C_\tau^4$ for each facet $\tau\in U$. Note
that apart from~$\tau$ and~$\tau'=\varphi'(\tau)$ there might be already a
triangulation induced on some cells
of~$\partial C_\tau^4$ via the triangulation of neighboring cells of~$C_\tau^4$. Fix a
5-coloring on the vertices of~$\{\tau,x_\tau\}\in N$, and color each vertex of~$\tau'$ with
the color of its preimage under~$\psi_{\{\tau,x_\tau\}}$.

\begin{lem}\label{lem:attaching1}
  The $5$-coloring of~$\tau$
  and~$\tau'$ can be extended to a $5$-coloring of the cells of~$\partial C_\tau^4$ already
  triangulated.
\end{lem}

\begin{proof}
  Let us call any strongly connected subcomplex of~$N$ with trivial
  group of projectivities which contains a facet $\s\in N$ a
  \emph{trivial domain} of~$\s$, and consider the trivial domain
  of~$\{\tau,x_\tau\}$
  \[
  O=\!\bigcup_{v\in\{\tau,x_\tau\}\setminus N_{\odd}}\!\Star_N(v),
  \]
  defined by the union of the stars of all
  vertices of $\{\tau,x_\tau\}$ not contained in~$N_{\odd}$. This
  is indeed a trivial domain if~$N$ is triangulated sufficiently fine
  (there are no identifications in~$\partial O$),
  since no star of an odd triangle is contained in~$O$, and since any facet
  path in~$O$ is contractable. For each cell~$C_f^{i+1}$ of~$\partial C_\tau^4$ already
  triangulated there is a facet~$\rho\in U$ in with $f=\tau\cap\rho$,
  and in the case $f\not\in N_{\odd}$ we have $\{\rho,x_\rho\}\in O$. Hence the
  5-coloring of~$\{\tau,x_\tau\}$ extends uniquely to the triangulation
  of~$C_f^{i+1}$. Furthermore, if there are two facets~$\rho$ and~$\overline{\rho}$
  with $f=\tau\cap\rho=\tau\cap\overline{\rho}$, both facets~$\rho$
  and~$\overline{\rho}$ produce the same coloring of the triangulation
  of~$C_f^{i+1}$ since $\Star_N(f)\subset O$, and since~$O$ is a trivial
  domain.

  In the case where $f\in N_{\odd}$ consider the subcomplex
  $\overline{O}=\bigcup_{v\in\{\tau,x_\tau\}}\Star_N(v)$, a regular
  neighborhood of~$\{\tau,x_\tau\}$. Assuming a sufficiently fine
  triangulation of~$N$ and that~$N_{\odd}$ is locally flat, we have
  $\overline{O}\cong\B^4$, $\overline{O}_{\odd}\cong\B^2$ with
  $\overline{O}_{\odd}\cap\partial\overline{O}\cong\Sph^1$, and
  $\Pi(\overline{O})\cong\Sigma_2$. Therefore 3 colors of the
  5-coloring of~$\{\tau,x_{\tau}\}$ corresponding to the three
  trivial orbits of $\Pi(\overline{O})$, and let us call these three
  colors the \emph{stable colors}.  Propagating the 5-coloring
  of~$\{\tau,x_\tau\}$ along any facets path in~$\overline{O}$
  from~$\{\tau,x_\tau\}$ to any facet $\{\rho,x_\rho\}\in
  \overline{O}$ with $f=\tau\cap\rho$ yields the same coloring for the
  triangulation of~$C_f^{i+1}$ using only the three stable colors,
  since the vertices of~$\{f,x_f\}$ correspond to trivial orbits of
  $\Pi(\overline{O})$.
\end{proof}

Now the partial triangulation and 5-coloring of~$\partial
C_\tau^4$ is extended to a triangulation and 5-coloring of the
entire cell~$C_\tau^4$ using Proposition~\ref{prop:extend_cw_complex}.
The triangulation and 5-coloring of~$C_\tau^4$ is extended in
two steps. First, let $f=\tau\cap N_{\odd}$, and
triangulate~$C_f^{i+1}$ applying
Proposition~\ref{prop:extend_cw_complex} using only the three stable
colors, unless, of course,~$C_f^{i+1}$ is already triangulated.
Then using Proposition~\ref{prop:extend_cw_complex} once more, the
triangulation and 5-coloring is extended to the entire
cell~$C_\tau^4$.

\begin{lem}\label{lem:attaching2}
  The odd subcomplex of $K\cup T\cup K'$ is
  $K_{\odd}\cup K'_{\odd}$, and
  the union of all triangles in $\bigcup_{f\in U\cap K_{\odd}}C_f^{i+1}$.
\end{lem}

\begin{proof}
  We first prove that a triangle~$t$ in the interior of a cell~$C_f^{i+1}$
  is even if $f\not\in K_{\odd}$. To this end let~$\tau\in U$ be a
  facet with $f\in \tau$ and let~$O$ be the trivial domain
  of~$\{\tau,x_\tau\}$ as described above. By construction of~$T$ there
  is a 5-coloring of the triangulation of $\bigcup_{\tau\in
    O}C_\tau^4\supset \Star_T(t)$, thus~$t$ is even.
  Any triangle~$t$ in~$U$, respectively~$U'$, is even in
  $K\cup T\cup K'$, since for any
  facet $\tau \in U$ the
  5-coloring of the cell~$C_\tau^4$ extends the 5-coloring
  of~$\{\tau,x_\tau\}$ and ~$\{\tau',x_{\tau'}\}$ by construction of~$T$, hence $\Star_{K\cup
    T\cup K'}(t)$ is 5-colorable and~$t$ is even.
  
  It remains to determine the parity of the triangles in
  the union $\bigcup_{f\in
    U\cap N_{\odd}}C_f^{i+1}$, which form a PL 2-manifold (with
  boundary) equivalent to
  $|U\cap N_{\odd}|\times[0,1]$, and we denote the triangles in question
  suggestively by~$T_\text{O}$. Let~$e$ be an interior
  co-dimension 3-face of a combinatorial manifold, hence we have
  $\Link(e)\cong\Sph^2$. It is immediate by double counting facet-ridge incidences
  in any simplicial pseudo manifold without boundary, that the number of facets is even,
  thus $\Link(e)$, and consequently $\Star(e)$ has an even number of
  facets.
  We double count the number of incidences of co-dimension
  2-faces~$\{e,x\}\in\Star(e)$ incident to~$e$, and facets
  of~$\Star(e)$
  \[
  \sum_{\{e,x\}\in\Star(e)}\hspace{-.2cm}\sharp\SetOf{\s\in\Star(e)}{\{e,x\}\subset\s}
  = \sum_{\s\in\Star(e)} 3.
  \]
  The left hand side equals the number of odd co-dimension 2-faces incident
  to~$e$ modulo~2, and the right hand side is even since there is an even number of facets
  $\s\in\Star(e)$.
  
  Returning to our triangulation $K\cup T\cup K'$, we have that
  any edge $e\not\in\partial T_\text{O}$ is contained in none or
  two odd triangles ($e$ is a ridge of the 2-manifold~$T_\text{O}$). Therefore if there is one odd
  triangle in a (strongly) connected
  component of~$T_\text{O}$, then all triangles in the connected
  component of~$T_\text{O}$ must be
  odd faces of~$T$, and each connected
  component of~$T_\text{O}$ intersects~$K_{\odd}$ in at least one 
  edge.
  Thus all triangles in
  $\bigcup_{f\in U\cap K_{\odd}}C_f^{i+1}$ are odd, and we proved $T_\text{O}=T_{\odd}$.
\end{proof}

\begin{lem}\label{lem:attaching3}
  The simplicial complex $K\cup T\cup K'$ is a combinatorial
  $4$-manifold. In particular, $K\cup T\cup K'$ is a nice simplicial complex. 
\end{lem}

\begin{proof}
  It suffices to prove that the vertex link of each vertex in~$T$ is
  a simplicial 3-sphere or simplicial 3-ball (in $K\cup T\cup K'$),
  and hence a combinatorial sphere, respectively ball. Let $f\in U$ be an $i$-face, $C^{i+1}_f$
  the corresponding closed cell of the regular cell decomposition
  of $|U|\times[0,1]$, and let~$v\in T$ be a vertex contained in the triangulation
  of $C^{i+1}_f$. Further let $g\in U$ be an $j$-face containing~$f$
  (thus $i\leq j$).
  In the case $v\in U\subset T$ (or $v\in U'\subset T$) we have
  \[
  D(v,g)=\big|\Link_T(v)\big| \cap C^{j+1}_g \cong \Cone\big(
  \partial \Star_U(f) \cap g\big),
  \]
  and otherwise
  \[
  D(v,g)=\big|\Link_T(v)\big| \cap C^{j+1}_g \cong \Susp \big(
  \partial \Star_U(f) \cap g\big).
  \]

  Observe that if $i=3$, that is,~$f$ is a ridge of~$N$ (and a facet of~$U$), then $\Link_T(v)$ is
  a 3-ball (if $v\in\partial T$) or 3-sphere completely contained in~$C^4_f$.
  Otherwise $D(v,g)$ is a 3-ball in the
  case $v\in U\cup U'$ as well as in the case $v\not \in U\cup U'$ for
  $j\not=0$. In the remaining case $v\not \in U\cup U'$ and~$j=0$ we have $D(v,g) \cong \Sph^0$. (Recall that $\Cone(\emptyset)\cong\B^0$ and
  $\Susp(\emptyset)\cong\Sph^0$ holds by definition.)

  For $i<3$ let $\tau,\tau'\in\Star_U(f)$ be facets intersecting in
  \mbox{$g=\tau\cap\tau'\supset f$}. Then the two 3-balls $D(v,\tau)$
  and $D(v,\tau')$ intersect in $D(v,g)$. Assume that
  $f\not\in\partial U$ holds. Since $\Star_U(f)$ is a combinatorial
  3-ball (and $\partial \Star_U(f)$ a combinatorial 2-sphere) we have
  \[
  \Link_T(v) \cong \bigcup_{\tau\in\Star_U(f)}D(v,\tau)\cong
  \Cone(\partial \Star_U(f))
  \]
  if $v\in U\cup U'$, and
  \[
  \Link_T(v) \cong \bigcup_{\tau\in\Star_U(f)}D(v,\tau) \cong
  \Susp(\partial \Star_U(f))\] otherwise.  The case $f\in\partial U$
  is treated similarly, except we consider the 2-ball
  $\closure(\partial \Star_U(f)\setminus\partial U)$ instead of the
  entire boundary of $\Star_U(f)$.  Thus~$T$ is a combinatorial
  4-manifold.

  It remains to prove that $\Link_{K\cap T \cap K'}(v)$ is a
  3-sphere or 3-ball for a vertex $v\in U\subset T$ (or
  $v\in U'\subset T$). This follows
  since $\Link_{K\cap T \cap K'}(v)$ is the union of the two
  combinatorial 3-balls $\Link_T(v)$ and $\Link_K(v)$,
  respectively $\Link_{K'}(v)$.
\end{proof}

Attaching~$\bigcup_m D_m$ to~$D_A$ producing $D_A\cup D_H$, and then
attaching~$D_B$ to~$D_A\cup D_H$ as described above completes the
construction of~$S$, a combinatorial manifold homeomorphic to~$\Sph^4$
(note the difference to a combinatorial 4-sphere). The partial
unfolding~$\widehat{S}$ is a combinatorial manifold by
Remark~\ref{rem:unf_commutes}.

It remains to verify conditions~(I) and~(II).  As for
condition~(I),~$T_{\odd}$ is homotopy equivalent to~$U_{\odd}$.
Further, any path around an odd triangle in the triangulation of some
cell~$C_f^{i+1}$, where~$f$ is an edge in $U\cap N_{\odd}$, is
homotopy equivalent to a path around the (unique) triangle
$\{f,x_f\}\in N_{\odd}$. This settles condition~(II).  We summarize
the construction in the following Theorem~\ref{thm:ch_4mf_main} which
states the main result of this paper.

\begin{thm}\label{thm:ch_4mf_main}
  For every closed oriented PL $4$-manifold~$M$ there is a combinatorial manifold~$S\cong\Sph^4$
  such that one of the connected components of the partial
  unfolding~$\widehat{S}$ of~$S$ is a combinatorial $4$-manifold PL-homeomorphic to~$M$. The
  projection $\widehat{S}\to S$ is a simple $4$-fold branched cover
  branched over a PL-surface with a finite number of cusp and node singularities.
\end{thm}



\bibliographystyle{mod_siam}
\bibliography{literature}

\end{document}